\documentclass[11pt]{article}
\usepackage{graphicx} 
\usepackage{amsmath,amssymb}
\usepackage{theorem}
\usepackage{epsfig}
\usepackage{verbatim}
\usepackage{graphicx}
\usepackage{subfigure}
\usepackage{cancel}
\usepackage{epstopdf}
\usepackage{a4wide}
\usepackage{pdflscape}

\usepackage{xurl}
\usepackage{hyperref}

\usepackage{float}

\usepackage{mathdots}
\usepackage{color}

\usepackage{tkz-euclide}

\usepackage[backend=biber,url=false,eprint=false,maxbibnames=99,giveninits=true,block=space, 
]{biblatex}
\newtheorem{theorem}{Theorem}[section]
\newtheorem{lemma}[theorem]{Lemma}

\newtheorem{prop}[theorem]{Proposition}
\newtheorem{corollary}[theorem]{Corollary}

\newtheorem{rem}[theorem]{Remark}

\newcommand{\RR}{\mathbb{R}}

\newcommand{\Z}{\mathbb{Z}}

\newcommand{\T}{\mathbb{T}}

\newcommand{\sgn}{{\rm sgn}}

\newcommand{\clau}{{\rm Cl}_2}

\newcommand{\dx}{\pa_x}

\newcommand{\al}{\alpha}

\newcommand{\ga}{\gamma}
\newcommand{\de}{\delta}
\newcommand{\ep}{\varepsilon}

\newcommand{\te}{\theta}

\newcommand{\Ical}{\mathcal{I}}

\newcommand{\Kcal}{\mathcal{K}}

\newcommand{\pa}{\partial}

\newcommand{\intpi}{\int_{-\pi}^\pi}

\newcommand{\norm}[1]{\Vert#1\Vert}
\newcommand{\ti}[1]{\widetilde{#1}}
\newcommand{\ha}[1]{\widehat{#1}}

\title{Stationary homogeneous solutions\\ for the inviscid SQG equation}
\author{Miguel M.G. Pascual-Caballo}
\date{}

\newcommand{\XC}[1]{C^{#1}(\T)}

\newcommand{\OpSUnfold}{\mathcal{S}}
\newcommand{\OpSGen}{S}
\newcommand{\OpS}{\mathcal{S}_n}

\newcommand{\OpNolinearGen}{A}

\newcommand{\OpSmooth}{\mathcal{R}_n}

\newcommand{\ConvSet}{V_M}
\newcommand{\ConvSubSet}{\ti{V}_M}

\newcommand{\kernel}{K_n}

\newenvironment{proofsth}[1]{\begin{trivlist} \item[] {\textbf{Proof of #1:}}}{\hfill $\Box$
                       \end{trivlist}} 

\newif\ifhideproofs

\ifhideproofs

\else
    \newenvironment{proof}{
    \begin{trivlist} \item[] {\bf Proof:}
    }{\hfill $\Box$\end{trivlist}}
\fi

\usepackage{lmodern}
\usepackage{titlesec}
\titleformat{\chapter}[block]
  {\normalfont\huge\bfseries}{\thechapter.}{0em}{\Huge}

\numberwithin{equation}{section}

\begin{document}

\maketitle
\setcounter{section}{0}
\begin{abstract}
    In this paper, we prove the existence of nontrivial stationary homogeneous solutions with infinite energy (unbounded at infinity) for the SQG equation. Our analysis also covers the existence of stationary solutions for the generalized De Gregorio equation $\partial_t w+\al u\,\dx w= w\,\dx u,\ \dx u=Hw$ with $\al>\tfrac{1}{2}$, where $H$ denotes the Hilbert transform.
\end{abstract}

\section{Introduction}
The inviscid Surface Quasi-geostrophic equation (SQG) in $\RR^2$ is given by the following transport problem
\begin{equation}\label{SQG}
    \left\{\begin{array}{l}
    \partial_t\theta+\nabla^\perp u\,\cdot\, \nabla\theta=0,
    \quad (x,t)\in\RR^2\times [0,\infty),\\
    u=-\Lambda^{-1}\theta,\\
    \theta(x,0)=\theta_0(x),
    \end{array}\right.
\end{equation}
where $\Lambda:=(-\Delta)^{\frac{1}{2}}$. No decay condition is imposed on $\theta$ at spatial infinity, since the homogeneous solutions considered in this paper are unbounded as $r\to\infty$. For such $\theta$, $\Lambda^{-1}$ is understood as in \cite{castro2020lifespan,elgindi2018symmetries}, in the sense made precise in Remark~\ref{RemDistSense} below.

This equation describes the temperature evolution in a general quasi-geostrophic framework, applicable to atmospheric and oceanic flows, assuming small Rossby and Ekman numbers, along with constant potential vorticity. We refer to \cite{CMT,Held_1995,MajdaBertozzi02,Pedlosky} for more details.

The mathematical study of this equation received a great boost with the works of Constantin, Majda, and Tabak, \cite{CMT} and \cite{CMT2}. In them, the authors revealed the similarities between SQG and 3D Euler, and proposed it as an interesting model in which to search for singularities.

There are several known results concerning the existence of solutions to the SQG model. In \cite{Resnick95}, Resnick proved global existence in $L^2$, and Marchand later extended this result in \cite{marchand} to $L^p$, $p>\tfrac{4}{3}$. Proofs for local well-posedness in $H^s,\ s>2$ can be found in \cite{dongho,CMT}. Wu also proved local well-posedness in \cite{Wu2005} for $C^\ga\cap L^p$, $\ga>1$ and $p>1$.
In \cite{CoIgNg,constantinnguyen2016,constantinnguyen2018}, the authors studied SQG on bounded domains.

Regarding non-uniqueness, Buckmaster, Shkoller, and Vicol \cite{BuckmasterVicol19} constructed weak solutions satisfying $\Lambda^{-1}\te\in C^\ga,\ \tfrac{1}{2}<\ga<\tfrac{4}{5}$. This result was extended to $\Lambda^{-1}\te\in C^\ga$ for every $\ga<1$ by Dai, Giri, and Radu \cite{DGR}, as well as by Isett and Looi \cite{IL}
(see also \cite{zhao2025onsage}).

In the forcing case, there are similar non-uniqueness results by Bulut, Huynh, and Palasek \cite{BKP23} and by Dai and Peng \cite{dai2023nonuniqueweak}. Furthermore, Castro, Faraco, Mengual, and Solera \cite{CFMS2024,CFMS2025} constructed non-unique solutions for $\theta,f\in L^1\cap L^p,\ 1<p<\infty$ and $\theta,f\in W^{s,p}\cap \dot{H}^{-\frac{1}{2}}$, $0\leq s < 1 + \tfrac{2}{p},\ 1\leq p\leq \infty$.

A central question in fluid mechanics is whether SQG develops finite-time singularities. Kiselev and Nazarov \cite{kiselevnazarov} showed the existence of arbitrarily small initial data in $H^s$, $s\geq 11$, whose corresponding solution's $H^s$-norm exceeds any prescribed bound within a finite time interval. Moreover,
Friedlander and Shvydkoy \cite{friedlandershvydkoy} found unstable eigenvalues for the linearization around a steady state.
Additionally, He and Kiselev \cite{hekiselev} showed that solutions exist for which
either the $C^2$-norm grows exponentially for all time or the solution blows up in
finite time. Several examples of global solutions with nontrivial dynamics have also been constructed. We refer to \cite{delPino,CCGSglobal,Godard,GS}.

C\'{o}rdoba and Mart\'{i}nez-Zoroa \cite{CMZsqg} showed
ill-posedness in $C^k$ for $k\geq 2$ and in $H^s$ with $\tfrac{3}{2} < s < 2$. See also the works of Jeong and Kim \cite{jeongkimsqg}, as well as Elgindi and Masmoudi \cite{tareknaderill}. In addition, C\'{o}rdoba, Mart\'{i}nez-Zoroa, and O\.{z}a\'{n}ski \cite{CMZOeuler} showed
that this loss of regularity is continuous in time, with the solution's Sobolev index
decreasing steadily from $s$ rather than dropping at a single instant.

In this paper, we will consider homogeneous solutions, which are of the form \begin{align}
\label{ansatz}
\theta(x_1,x_2,t)=rf(x,t),
\end{align}
where $ (r,x)\in\RR^+\times \T$ are the usual polar coordinates, and $\T:=\RR/2\pi\Z$ denotes the one-dimensional torus, identified with $(-\pi,\pi]$.

Introducing the ansatz \eqref{ansatz} into \eqref{SQG}, one obtains the following equation for a function $f:\T\times[0,\infty)\to\RR$ whose Fourier coefficients satisfy $\ha{f}_k=0$ for $|k|\leq 2$,
\begin{equation}\label{EqnCircleDynamic}
\partial_t f+2\OpSUnfold f\, \partial_xf-\partial_x\OpSUnfold f\,f=0,
\end{equation}
where $\OpSUnfold$ is a Fourier multiplier defined by
\begin{align}\label{EqnDefUnfolded}
    (\ha{\OpSUnfold f})_k:=\ha{\OpSUnfold}_k\ha{f}_k:=-
        \frac{1}{|k|}\frac{k^2-1}{k^2-4}
        \ha{f}_k,\quad |k|\geq3.
\end{align}
Detailed computations leading to this expression can be found in \cite{castro2020lifespan} and \cite{elgindi2018symmetries}.

The purpose of this paper is to prove the existence of nontrivial stationary solutions to \eqref{EqnCircleDynamic}, which is a first step toward studying the long-time dynamics of nearby solutions or constructing non-uniqueness examples.

\begin{theorem}\label{MainTheorem}
    For each $n\geq 3$, there exists an odd $\frac{2\pi}{n}$-periodic function $f$, belonging to $\XC{\frac{1}{2}}$, such that $\ha{f}_n\neq0$, solving \eqref{EqnCircleDynamic}. In addition, the $n$-fold symmetric function $\Theta:\RR^2\to\RR$ given by
    \begin{align*}
        \Theta(x_1,x_2):=rf(x)
    \end{align*}
    is a stationary solution of the SQG equation~\eqref{SQG}. Moreover, if $q\in L^2(\T)$ is an odd $\frac{2\pi}{n}$-periodic solution of \eqref{EqnCircleDynamic}, and $\OpSUnfold q$ does not vanish on $\big(0,\tfrac{\pi}{n}\big)$, then $q$ is a multiple of $f$.
\end{theorem}

\begin{rem}\label{RemDistSense}
    The solutions constructed here do not depend on time, so \eqref{SQG} reduces to
    \begin{align*}
        \nabla^\perp u\cdot\nabla\Theta=0, \qquad u=-\Lambda^{-1}\Theta.
    \end{align*}
    This is understood in the distributional sense, that is,
    \begin{align*}
        \int_{\RR^2} \Theta\,\nabla^\perp u\cdot\nabla\psi\,dx=0
        \qquad\text{for every }\psi\in C_c^\infty(\RR^2),
    \end{align*}
    which comes from $\nabla^\perp u\cdot\nabla\Theta=\nabla\cdot\big(\Theta\nabla^\perp u\big)$, valid because $\nabla^\perp u$ is divergence free. Since $\nabla\psi$ is bounded and supported in a compact set, the integral is absolutely convergent as soon as $\Theta\nabla^\perp u$ is locally integrable.
\end{rem}

An earlier version of this work contained the existence result only. Once it was completed, we
became aware of the paper \cite{SteadySerrano} by Abe, Gómez-Serrano, and Jeong, and we agreed to
synchronize its upload with theirs on arXiv. Their results for SQG also included solutions of the
form $\theta=rf(x)$, obtained through variational methods, whereas our proof relied on a
fixed-point argument (see Proposition~\ref{PropGeneralExistence}). Our results additionally covered the generalized De Gregorio equation
\eqref{DeGregGen} below on $\T$, not addressed in \cite{SteadySerrano}.

In the process of proving Theorem~\ref{MainTheorem}, we prove a more general result. Specifically, we consider the equation
\begin{equation}\label{EqnStatFolded}
    \al \OpSUnfold f\, \partial_x f=\partial_x\OpSUnfold f\,f,
\end{equation}
with $\al>\tfrac{1}{2}$. We obtain nontrivial, odd, and $\frac{2\pi}{n}$-periodic solutions for each $n\geq 3$.

The family of equations
\begin{equation}\label{EqnCircleDynamica}
\partial_t f+\al \OpSUnfold f\, \partial_xf-\partial_x\OpSUnfold f\,f=0,
\end{equation}
with $\al\in \RR$, is intimately related to the generalized De Gregorio equation,
\begin{equation}\label{DeGregGen}
    \left\{\begin{array}{l}
    \partial_t w+\al u\,\dx w= w\,\dx u,\quad (x,t)\in\RR\times [0,\infty),\ \al\in\RR,\\
    u=-\Lambda^{-1} w,\\
     w(x,0)= w_0(x).
    \end{array}\right.
\end{equation}
Here $\Lambda:=(-\Delta)^{\frac{1}{2}}$. Equation \eqref{DeGregGen} with $\al=1$ was proposed by De Gregorio \cite{DeGregorio1996,DeGregorio1990} as a
one-dimensional model that captures the
competition between advection and vortex stretching in the Euler equation. Okamoto, Sakajo, and Wunsch \cite{Okamoto2008} introduced the parameter $\al\in\RR$ into the De Gregorio equation, obtaining the generalized model, also known in the literature as the generalized Constantin--Lax--Majda (gCLM) equation.

If we fix $\al=-1$ in \eqref{DeGregGen}, we obtain the C\'{o}rdoba--C\'{o}rdoba--Fontelos model \cite{cordoba2007formation}. Furthermore, the case $\al=0$ corresponds to the Constantin--Lax--Majda equation, see \cite{CLMModel}.

Equations \eqref{EqnCircleDynamica} and \eqref{DeGregGen} are closely related because
\begin{align*}
     \OpSUnfold=-\Lambda^{-1}+\text{smoothing operator}.
\end{align*}

Okamoto, Sakajo, and Wunsch \cite{OkaSaWu} studied the existence of nontrivial steady states of $\eqref{DeGregGen}$ in $C^{\ga}(\T)$, with $0<\al\leq1$. In this paper, we prove the existence of nontrivial solutions to \eqref{EqnStatFolded} for every $\al>\tfrac{1}{2}$, in $C^{\frac{1}{\al}}(\T)$ when $\al>1$, and in $C^{1,\frac{1}{\al}-1}(\T)$ when $\al<1$.
Our argument, inspired by \cite{OkaSaWu}, is also valid for the generalized De Gregorio equation $\eqref{DeGregGen}$.

\subsection*{Note added}
The version of this work posted simultaneously with \cite{SteadySerrano} contained the
existence result only. The argument has since been modified so that it also yields the
uniqueness statement in Theorem~\ref{MainTheorem}.

\section{Construction of stationary solutions}

Here and below, $\dx^{-1}$ denotes the zero-mean primitive on $\T$ and $H$ the Hilbert transform. On mean-zero functions they are given by the Fourier multipliers
\begin{align*}
    (\ha{\dx^{-1}f})_k:=\frac{1}{ik}\ha{f}_k
    \qquad
    \text{and}
    \qquad(\ha{Hf})_k:=-i\,\sgn(k)\,\ha{f}_k\qquad\text{for } k\neq0.
\end{align*}
Notice that the operator $\Lambda^{-1}$ in \eqref{DeGregGen} can be written as $\Lambda^{-1}=-\dx^{-1}\!H=-H\dx^{-1}$.
We also write $A\lesssim B$ to mean $A\leq CB$ for a constant $C$ that may change from one line to the next.

We seek nontrivial solutions to \eqref{EqnStatFolded}. We will assume that $f\in L^2(\T)$ is odd and $\frac{2\pi}{n}$-periodic, $n\geq 3$, so we can expand

\begin{align*}
    f(x)
    =
    \sum_{k\geq 1}
        \ha{f}^{\rm s}_{nk}
        \sin(nkx),
\end{align*}
where the sine coefficients of $f$ are given by
\begin{align*}
    \ha{f}^{\rm s}_k:=\frac{1}{\pi}\int_{-\pi}^{\pi}f(x)\sin(kx)\,dx=\frac{2}{\pi}\int_0^\pi f(x)\sin(kx)\,dx.
\end{align*}

\begin{rem}\label{RemRescaling}
    Suppose that there is an odd stationary solution of \eqref{DeGregGen}, with nonzero first Fourier coefficient. From such a solution one can construct new ones by rescaling the period.
    This follows from the fact that $-\Lambda^{-1}=\dx^{-1}\!H$ and that, for any mean-zero function $g\in C(\T)$, if we set $g_n:=g(nx)$, then
    \begin{align*}
        \dx^{-1}\!Hg_n(x)=\frac{1}{n}\dx^{-1}\!Hg(nx),\quad n\geq 1.
    \end{align*}
    Notice that this is not the case for $\OpSUnfold$, since $\OpSUnfold g_n(x)\neq \frac{1}{n}\OpSUnfold g(nx)$, and this is what makes our analysis different. We are nevertheless able to prove the theorem for every $n\geq 3$, thereby generating this family of solutions.
\end{rem}

For $n\geq 3$, we denote by $\OpS$ the Fourier multiplier given by
\begin{equation}
\label{EqnStatUnFolded}
    (\ha{\OpS w})_k
    :=(\ha{\OpS})_k\ha{w}_k
    :=
    -\frac{1}{|k|}\frac{n^2k^2-1}{n^2k^2-4}
    \ha{w}_k,
    \qquad k\neq 0.
\end{equation}
The following lemma collects the properties of $\OpS$ and $\OpSUnfold$ that we will use, and the relationship between them.
\begin{lemma}
\label{Lemma_FourierUnfoldFold}
    Let $n\geq 3$. Both $\OpS$ and $\OpSUnfold$ take odd functions to odd functions. Moreover, if $f\in L^2(\T)$ is odd and $\frac{2\pi}{n}$-periodic, then, for $g(x):=f\big(\tfrac{x}{n}\big)$, both $\OpS g$ and $\OpSUnfold f$ belong to $H^1(\T)\subset C(\T)$, and
    \begin{align*}
        \OpS g(x)=n\,\OpSUnfold f\big(\tfrac{x}{n}\big)
        \qquad\text{for every }x\in\T.
    \end{align*}
\end{lemma}
\begin{proof}
    If $w$ is odd then $\ha{w}_{-k}=-\ha{w}_k$, and since $(\ha{\OpS})_k$ depends on $k$ only through $|k|$,
    \begin{align*}
        (\ha{\OpS w})_{-k}=(\ha{\OpS})_{-k}\ha{w}_{-k}=-(\ha{\OpS})_k\ha{w}_k=-(\ha{\OpS w})_k,
    \end{align*}
    so $\OpS w$ is odd. The same applies to $\OpSUnfold$, whose symbol is defined for $|k|\geq 3$.

    Let now $f\in L^2(\T)$ be odd and $\frac{2\pi}{n}$-periodic. It can be written as the sine Fourier series
    \begin{align*}
        f(x)=\sum_{k\geq 1}\ha{f}^{\rm s}_{nk}\sin(knx),
    \end{align*}
    converging in $L^2(\T)$, and hence
    \begin{align}\label{EqnExpansiong}
        g(x):=f\big(\tfrac{x}{n}\big)=\sum_{k\geq 1}\ha{f}^{\rm s}_{nk}\sin(kx)
    \end{align}
    in $L^2(\T)$ as well.

    Both symbols are defined and bounded, see \eqref{EqnDefUnfolded} and \eqref{EqnStatUnFolded}, so that $\OpS g$ and $\OpSUnfold f$ are the $L^2(\T)$ limits of the corresponding partial sums.
    Now, observe that
    \begin{align}\label{FourierUnfoldFold}
        \ha{\OpSUnfold}_{kn}=
        -\frac{1}{n|k|}\frac{n^2k^2-1}{n^2k^2-4}=
        \frac{1}{n}(\ha{\OpS})_k.
    \end{align}
    Since $f$ is $\frac{2\pi}{n}$-periodic, its only nonzero Fourier coefficients occur at multiples of $n$, so, using \eqref{EqnExpansiong} and \eqref{FourierUnfoldFold},
    \begin{align*}
        \OpS g(x)=\sum_{k\geq 1}(\ha{\OpS})_k\ha{g}^{\rm s}_k\sin(kx)=\sum_{k\geq 1}(\ha{\OpS})_k\ha{f}^{\rm s}_{nk}\sin(kx)=n\sum_{k\geq 1}\ha{\OpSUnfold}_{kn}\ha{f}^{\rm s}_{nk}\sin(kx)=n\,\OpSUnfold f\big(\tfrac{x}{n}\big)
    \end{align*}
    as an identity in $L^2(\T)$. Moreover, $(\ha{\OpS})_k=O\big(\tfrac{1}{k}\big)$ and $\ha{\OpSUnfold}_k=O\big(\tfrac{1}{k}\big)$, so both sides above in fact belong to $H^1(\T)\subset C(\T)$, and the identity holds for every $x\in\T$.
\end{proof}

Then, from Lemma~\ref{Lemma_FourierUnfoldFold}, equation \eqref{EqnStatFolded} becomes
\begin{align}
\label{EqnStatFoldedg}
    \al \OpS g\, \dx g=g\, \dx \OpS g,
\end{align}
for a $2\pi$-periodic, mean zero and odd function $g\in L^2(\T)$.

As $\OpS$ maps sine series to sine series, it is natural to let it act on functions on $[0,\pi]$, by applying it to the odd extension and restricting the result back to $[0,\pi]$.
Following the strategy outlined in \cite{OkaSaWu}, instead of \eqref{EqnStatFoldedg}, we will work with the equation
\begin{equation}
\label{EqnIntegral}
    v(x) = -\OpS\big(v^{\frac{1}{\al}}\big)(x)\quad\text{for every } x\in[0,\pi],
\end{equation}
where $v$ is a nonnegative function. Once we have a solution $v$ of \eqref{EqnIntegral}, we also have a solution $g$ of \eqref{EqnStatFoldedg}, as the following lemma states.
\begin{lemma}
\label{FromFixedPointToEQN}
    Let $\al>\tfrac{1}{2}$ with $\al\neq 1$, let $n\geq 3$, and let $v\in C^{1,\ga}([0,\pi]),\ \ga\in(0,1)$ with $\ga\geq\max\big\{0,\tfrac{1}{\al}-1\big\}$, satisfying $v(x)>0$ for every $x\in(0,\pi)$ and $v(0)=v(\pi)=0$, be a solution of \eqref{EqnIntegral}. Let $g\in C^\ga([0,\pi])$ be such that $\OpS g=v$. Then, the odd extension of $g$ is in $C^{\frac{1}{\al}}(\T)$ if $\al>1$, and in $C^{1,\frac{1}{\al}-1}(\T)$ if $\al<1$, and in both cases solves \eqref{EqnStatFoldedg}.
\end{lemma}

\begin{rem}\label{RemWeakSense}
    When we say that $g$ solves \eqref{EqnStatFoldedg}, it is understood in the divergence form
    \begin{align*}
        \al \dx\big(g\,\OpS g\big)=(\al+1)\,g\,\dx\OpS g,
    \end{align*}
    in which $g$ is never differentiated, that is,
    \begin{align*}
        -\al\intpi g(x)\,\OpS g(x)\,\dx\psi(x) \,dx
        =(\al+1)\intpi g(x)\,\dx\OpS g(x)\,\psi(x) \,dx
        \quad\text{for every }\psi\in C^\infty(\T).
    \end{align*}
    Since $(\ha{\OpS})_k=O\big(\tfrac{1}{k}\big)$, the operator $\OpS$ maps $L^2(\T)$ into $H^1(\T)\subset L^2(\T)$, so $\OpS g,\dx\OpS g,g\in L^2(\T)$. By Cauchy-Schwarz, both $g\,\OpS g$ and $g\,\dx\OpS g$, each a product of two $L^2(\T)$ factors, belong to $L^1(\T)$ and the two integrals converge. Equation \eqref{EqnStatFolded} is understood in the same sense.
\end{rem}

\begin{proof}
    Both $g$ and $v^{\frac{1}{\al}}$ are continuous on $[0,\pi]$, so their odd extensions belong to $L^2(\T)$ and are given by sine series. Since $(\ha{\OpS})_k\neq 0$ for every $k\geq 1$, the operator $\OpS$ is injective on sine series, so \eqref{EqnIntegral} gives
    \begin{align}\label{gisminusvpower}
        g=-v^{\frac{1}{\al}}\quad\text{on }[0,\pi].
    \end{align}

    Suppose first that $\al>1$, so that $\tfrac{1}{\al}\in(0,1)$. For $a\geq b\geq0$ and $t\in(0,1]$, the function $\phi(x):=(b+x)^t-x^t-b^t$ satisfies $\phi(0)=0$ and, for $x>0$,
    \begin{align*}
        \pa_x\phi(x)=t\big[(b+x)^{t-1}-x^{t-1}\big]\leq0,
    \end{align*}
    since $t-1\leq0$ and $b+x\geq x$, so $\phi$ is non-increasing and $\phi(a-b)\leq\phi(0)=0$, that is, $a^t-b^t\leq(a-b)^t$, and the case $b>a$ follows by symmetry, so
    \begin{align}\label{IneqSubadditivePower}
        \big|a^t-b^t\big|\leq|a-b|^t
        \qquad\text{for every }a,b\geq 0\text{ and every }t\in(0,1].
    \end{align}
    Since $v\in C^{1,\ga}([0,\pi])$, we obtain from \eqref{gisminusvpower} and \eqref{IneqSubadditivePower}, for every $x,y\in[0,\pi]$,
    \begin{align*}
        |g(x)-g(y)|
        =
        \big|v(x)^{\frac{1}{\al}}-v(y)^{\frac{1}{\al}}\big|
        \leq
        |v(x)-v(y)|^{\frac{1}{\al}}
        \leq
        \norm{\dx v}_{C([0,\pi])}^{\frac{1}{\al}}|x-y|^{\frac{1}{\al}},
    \end{align*}
    so $g\in C^{\frac{1}{\al}}([0,\pi])$.

    Since $g=-v^{\frac{1}{\al}}$ by \eqref{gisminusvpower} and $v(0)=v(\pi)=0$, we have $g(0) = g(\pi) = 0$, so the odd extension $\ti{g}$ of $g$ to the torus vanishes at $0$ and at $\pi$. On $[0,\pi]$, $\ti g=g\in C^{\frac{1}{\al}}$ as shown above, and on $[-\pi,0]$, $\ti g(x)=-g(-x)$ is Hölder with the same constant. It is enough to bound the distance between points separated by $0$ or by $\pi$. For $y < 0 < x$,
    \begin{align*}
        |\ti{g}(x) - \ti{g}(y)|
        &\leq |\ti{g}(x) - \ti{g}(0)| + |\ti{g}(0) - \ti{g}(y)| =
        |\ti{g}(x)| + |\ti{g}(y)| \\
        &\lesssim |x|^{\frac{1}{\al}} + |y|^{\frac{1}{\al}}
        \leq 2 |x-y|^{\frac{1}{\al}}.
    \end{align*}
    Here, in the last step, we have used that $|x-y| \geq |x|, |y|$. The same argument applies at the crossing at $\pi$. This gives the statement for $\al>1$.

    Suppose now that $\al\in\big(\tfrac{1}{2},1\big)$, so that $\tfrac{1}{\al}-1\in(0,1)$. By \eqref{gisminusvpower} and the chain rule,
    \begin{align*}
        \dx g=-\frac{1}{\al}\,v^{\frac{1}{\al}-1}\,\dx v,
    \end{align*}
    which is continuous on $[0,\pi]$ as a product of two continuous functions and vanishes at $0,\pi$. Applying \eqref{IneqSubadditivePower} with $t=\tfrac{1}{\al}-1\in(0,1)$, we obtain
    \begin{align*}
        \big|v(x)^{\frac{1}{\al}-1}-v(y)^{\frac{1}{\al}-1}\big|
        &\leq
        |v(x)-v(y)|^{\frac{1}{\al}-1}
        \leq
        \norm{\dx v}_{C([0,\pi])}^{\frac{1}{\al}-1}|x-y|^{\frac{1}{\al}-1},
    \end{align*}
    where the last inequality is the mean value theorem. Hence $v^{\frac{1}{\al}-1}\in C^{\frac{1}{\al}-1}([0,\pi])$.

    Both $v^{\frac{1}{\al}-1}$ and $\dx v$ are bounded on $[0,\pi]$, and
    \begin{align*}
        \big|\dx v(x)-\dx v(y)\big|
        &\leq
        \norm{\dx v}_{C^{\ga}([0,\pi])}|x-y|^{\ga}
        \leq
        \norm{\dx v}_{C^{\ga}([0,\pi])}\,\pi^{\ga-\frac{1}{\al}+1}|x-y|^{\frac{1}{\al}-1},
    \end{align*}
    so $\dx v\in C^{\frac{1}{\al}-1}([0,\pi])$, and the decomposition
    \begin{align*}
        \big|\dx g(x)-\dx g(y)\big|
        &\leq
        \tfrac{1}{\al}\,v(x)^{\frac{1}{\al}-1}\big|\dx v(x)-\dx v(y)\big|
        +
        \tfrac{1}{\al}\,|\dx v(y)|\,\big|v(x)^{\frac{1}{\al}-1}-v(y)^{\frac{1}{\al}-1}\big|
        \\
        &\lesssim
        |x-y|^{\frac{1}{\al}-1}
    \end{align*}
    gives $\dx g\in C^{\frac{1}{\al}-1}([0,\pi])$, that is, $g\in C^{1,\frac{1}{\al}-1}([0,\pi])$.

    For the odd extension, $g(0)=g(\pi)=0$ implies that $\ti g$ is continuous, and $\dx\ti g$ is even about $0$ and about $\pi$, hence continuous on $\T$. For $y<0<x$, $\dx\ti g(x)=\dx g(x)$ and $\dx\ti g(y)=\dx g(-y)$, so
    \begin{align*}
        \big|\dx\ti g(x)-\dx\ti g(y)\big|
        &=
        \big|\dx g(x)-\dx g(-y)\big|\\
        &\leq
        \norm{\dx g}_{C^{\frac{1}{\al}-1}([0,\pi])}|x+y|^{\frac{1}{\al}-1}
        \leq
        \norm{\dx g}_{C^{\frac{1}{\al}-1}([0,\pi])}|x-y|^{\frac{1}{\al}-1}.
    \end{align*}
    Since the same argument applies at $x=\pi$, we have $\ti g\in C^{1,\frac{1}{\al}-1}(\T)$.

    We now prove that $\ti g$ solves \eqref{EqnStatFoldedg} on $\T$, in the sense of Remark~\ref{RemWeakSense}.

    From \eqref{gisminusvpower} and $\OpS g=v$, we obtain
    \begin{equation}
        |g|^\al = |\OpS g|.
        \label{identity-g}
    \end{equation}
    Since $\OpS g=v>0$ on $(0,\pi)$, and both $g$ and $\OpS g$ are odd, it follows from \eqref{identity-g} that $g\neq 0$ and $\OpS g\neq 0$ on $\T-\{0,\pi\}$. Differentiating \eqref{identity-g} gives
    \begin{equation*}
        \al |g|^{\al-1}\sgn(g)\,\dx g
        =
        \sgn(\OpS g)\,\dx\OpS g.
    \end{equation*}
    Multiplying by $g$, we find
    \begin{equation}
        \al |g|^\al \dx g
        =
        g\,\sgn(\OpS g)\,\dx\OpS g.
        \label{chainImplications}
    \end{equation}

    Then, using $|g|^\al=|\OpS g|$ in
    \eqref{chainImplications}, we get
    \begin{align*}
        \al |\OpS g|\,\dx g
        =
        g\,\sgn(\OpS g)\,\dx\OpS g.
    \end{align*}
    Multiplying by $\sgn(\OpS g)$, this yields, on $\T-\{0,\pi\}$,
    \begin{equation}
        \al \OpS g\,\dx g
        =
        g\,\dx\OpS g.
        \label{IdentityPointwise}
    \end{equation}

    It remains to control the two excluded points, $0$ and $\pi$. Since $\OpS g=v\in C^{1,\ga}([0,\pi])$ vanishes at $0$ and at $\pi$, $\ti{\OpS g}\in C^{1,\ga}(\T)$.

    For $\al<1$ nothing more is needed, since $\ti g\in C^{1,\frac{1}{\al}-1}(\T)$ too, so both sides of \eqref{IdentityPointwise} are continuous on $\T$. For $\al>1$, note that $\ti g\,\ti{\OpS g}$ is continuous on $\T$ and that, away from $0$ and $\pi$, \eqref{IdentityPointwise} gives
    \begin{align*}
        \al \dx\big(\ti g\,\ti{\OpS g}\big)=(\al+1)\,\ti g\,\dx\ti{\OpS g},
    \end{align*}
    classically on $(0,\pi)$ and on $(-\pi,0)$. Since $\dx\big(g\,\OpS g\big)=\big(1+\tfrac1\al\big)g\,\dx\OpS g$ is bounded on $(0,\pi)$ and on $(-\pi,0)$, $\ti g\,\ti{\OpS g}$ is Lipschitz on each closed arc, hence absolutely continuous, so the identity extends across $0$ and $\pi$ and integration by parts against $\psi\in C^\infty(\T)$ is valid on all of $\T$, giving exactly the divergence form of \eqref{EqnStatFoldedg} in the sense of Remark~\ref{RemWeakSense}.
\end{proof}

\subsection{Schauder's fixed point argument}
We adapt the existence argument of \cite{OkaSaWu} to our setting.
The strategy is still based on Schauder's fixed point theorem \cite{Schauder1930},
applied to the set of functions normalized so that $\int_0^\pi v(x)\sin x\,dx=1$.
We first state the result in a general form, valid for other operators. It then suffices to check that
$-\OpS$ satisfies the required assumptions.
\begin{prop}\label{PropGeneralExistence}
    Let $M>0$ be such that the set $\ConvSet$ given by
        \begin{align*}
        \ConvSet:=\left\{
        v \in C([0,\pi]):\
        0\leq v\leq M,\
        \int_0^\pi v(x) \sin x \,dx = 1
        \right\},
    \end{align*}
    is non-empty.

    Let $\al>0$, with $\al \neq 1$, and let
    \begin{align*}
        \OpSGen:C([0,\pi])\to C([0,\pi])
    \end{align*}
    be a linear compact operator. Define
    \begin{align*}
        \OpNolinearGen v
        :=
        \OpSGen\big(v^\frac{1}{\al}\big).
    \end{align*}
    Assume that $\OpNolinearGen$ satisfies the following properties
    \begin{enumerate}
        \item\label{enume_first} There exists $\rho>0$ such that, for every $v\in\ConvSet$,
        \begin{align*}
            \int_0^\pi \OpNolinearGen v(x)\sin x\,dx
        =
        \rho\int_0^\pi v(x)^{\frac{1}{\al}}\sin x\,dx.
        \end{align*}
        \item\label{enume_second} For every $v\in\ConvSet$ and every $x\in[0,\pi]$,
        \begin{align*}
            0\leq
        \OpNolinearGen v(x)\leq
        M\int_0^\pi \OpNolinearGen v(y) \sin y \,dy.
        \end{align*}
    \end{enumerate}
    Then, $\OpNolinearGen$ has a nonzero fixed point in $C([0,\pi])$, which is a multiple of an element of $\ConvSet$.
\end{prop}
\begin{rem}\label{RemWhyM}
    To ensure that $\ConvSet$ is bounded, we have introduced the constant $M$, which does not appear in \cite{OkaSaWu}. This seems to be necessary to apply Schauder's fixed point theorem.
\end{rem}
\begin{proof}
    Let $F$ be the operator given by
    \begin{align*}
        Fv:=\frac{\OpNolinearGen v}{\int_0^\pi \OpNolinearGen v(x) \sin x \,dx}.
    \end{align*}
    Now, we observe that $\ConvSet$ is a non-empty, closed, bounded, and convex subset of $C([0,\pi])$. We also note that $Fv\in\ConvSet$ since $0\leq Fv\leq M$ and
    \begin{align*}
        \int_0^\pi Fv(x) \sin x \,dx=1.
    \end{align*}

    To apply Schauder's fixed point theorem, we need to show that $F$ is compact. Notice that $\OpSGen$ is compact, so is $\OpNolinearGen$. Therefore, it remains to check that
    \begin{align*}
        \inf\limits_{v\in\ConvSet}\int_0^\pi \OpNolinearGen v(x) \sin x \,dx=
        \rho\inf\limits_{v\in\ConvSet}
        \int_0^\pi
            v(x)^{\frac{1}{\al}}
            \sin x
        \,dx>0.
    \end{align*}

    For $\al<1$ we can use the same Hölder estimate as in \cite{OkaSaWu}:
    \begin{align}
        1=\int_0^\pi
        v(x) \sin x \,dx
        &\leq
        \underbrace{\left(\int_0^{\pi}\sin x \,dx\right)^{1-\al}}_{=\,2^{1-\al}}
        \left(\int_0^{\pi}v(x)^{\frac{1}{\al}} \sin x \,dx\right)^{\al}.
        \label{IneqDevelLowerBoundAlSma}
    \end{align}

    For the case $\alpha > 1$, we note that $\tfrac{1}{\al}-1<0$, so
    \begin{align*}
        v^{\frac{1}{\al}-1}\geq M^{\frac{1}{\al}-1}.
    \end{align*}
    Hence,
    \begin{align*}
        \int_0^{\pi}v(x)^{\frac{1}{\al}} \sin x\,dx\geq
        M^{\frac{1}{\al}-1}\int_0^{\pi}v(x) \sin x\,dx=M^{\frac{1}{\al}-1}.
    \end{align*}

    By Schauder's fixed point theorem, there exists $v \in \ConvSet$ which is a fixed point of $F$.
    Since $\alpha \neq 1$, we can construct a nonzero multiple $\Gamma v$ so that $\Gamma v$ is a fixed point of $\OpNolinearGen$. Indeed, set $c:=\int_0^\pi \OpNolinearGen v(x) \sin x\,dx$, so that the identity $Fv=v$ reads $\OpNolinearGen v=c\,v$. Since $\OpNolinearGen(\Gamma v)=\Gamma^{\frac{1}{\al}}\OpNolinearGen v$, the choice $\Gamma:=c^{\frac{\al}{\al-1}}$ gives $\OpNolinearGen(\Gamma v)=\Gamma v$.
\end{proof}

\begin{rem}
The proposition applies for every $\al>0$ with $\al\neq1$, but in this paper we verify its
hypotheses only for $\al\in\big(\tfrac12,1\big)$ (see Section~\ref{checking}). The case
$\al>1$ is instead treated directly via the contraction argument of
Proposition~\ref{PropGeneralUniqueness}.
\end{rem}

\begin{rem}\label{RemOddExtension}
    In Section~\ref{checking} we will take $\OpSGen=-\OpS$. For $w\in C([0,\pi])$, the function $\OpS w$ is the restriction to $[0,\pi]$ of an odd function of $H^1(\T)\subset C(\T)$, since $(\ha{\OpS})_k=O\big(\tfrac{1}{|k|}\big)$, so it vanishes at $0$ and at $\pi$. Any nonnegative fixed point $v\in C([0,\pi])$ of $\OpNolinearGen$ satisfies
    \begin{align*}
        v=\OpSGen\big(v^{\frac{1}{\al}}\big),
    \end{align*}
    and since $v$ is continuous and nonnegative, $v^{\frac{1}{\al}}\in C([0,\pi])$ too, so $v$ is likewise the restriction of an odd $H^1(\T)$ function, hence vanishes at $0$ and $\pi$ and its odd extension is continuous on $\T$.
\end{rem}

\begin{rem}
Finally, we observe that if there exists a non-empty, closed, and convex subset $\ConvSubSet \subset \ConvSet$ such that $F(\ConvSubSet) \subset \ConvSubSet$,
then Schauder's theorem can be applied to this smaller set, thereby yielding additional properties of the fixed point.
For instance, in our case we can intersect $\ConvSet$ with the space of functions satisfying
\begin{align*}
    v\!\left(\tfrac{\pi}{2}+x\right)=v\!\left(\tfrac{\pi}{2}-x\right).
\end{align*}
This subset is preserved by $F$ whenever $\OpSGen$ commutes with the reflection $Rv(x):=v(\pi-x)$, since then
\begin{align*}
    \OpNolinearGen Rv=\OpSGen\big((Rv)^{\frac{1}{\al}}\big)=\OpSGen R\big(v^{\frac{1}{\al}}\big)=R\OpSGen\big(v^{\frac{1}{\al}}\big)=R\OpNolinearGen v,
\end{align*}
and, since $\sin(\pi-y)=\sin y$, the normalizing integral $\int_0^\pi\OpNolinearGen v(y)\sin y\,dy$ is unchanged as well, so $F(Rv)=R(Fv)$.
Hence, $v$ admits a Fourier representation of the form
\begin{align*}
    v(x)=\sum_{k\geq 0}\ha{v}^{\rm s}_{2k+1}\sin((2k+1)x).
\end{align*}
\end{rem}

\subsection{Banach's fixed point argument}\label{uniqueness}

Proposition~\ref{PropGeneralExistence} says nothing about how many fixed points $\OpNolinearGen$ has. For $\al>1$ the exponent $\tfrac{1}{\al}$ is smaller than $1$, and this makes $\OpNolinearGen$ contract distances between functions that are comparable to the sine function, so we can run a fixed point argument as in \cite{Bushell1973,Thompson1963} instead. We write
\begin{align}\label{DefComparableSet}
    \Kcal:=\left\{v\in C([0,\pi]):\ \xi\sin x\leq v(x)\leq\Xi\sin x\ \text{ on }[0,\pi]\ \text{ for some }0<\xi\leq\Xi\right\},
\end{align}
and, for $v,w\in\Kcal$,
\begin{align}\label{DefComparableDistance}
    d(v,w):=\sup_{x\in(0,\pi)}\left|\log\frac{v(x)}{w(x)}\right|.
\end{align}

\begin{lemma}\label{LemmaComparableComplete}
    $(\Kcal,d)$ is a complete metric space.
\end{lemma}
\begin{proof}
    Let $C_b((0,\pi))$ denote the space of bounded continuous functions on $(0,\pi)$, with the supremum norm. Notice that the map
    \begin{align}\label{MapLogSin}
        v\mapsto\log\!\left(\frac{v}{\sin}\right)
    \end{align}
    takes $\Kcal$ into $C_b((0,\pi))$. Moreover, it preserves distances, since
    \begin{align*}
        \left\|\log\!\left(\frac{v}{\sin}\right)-\log\!\left(\frac{w}{\sin}\right)\right\|_{C_b((0,\pi))}=\sup_{x\in(0,\pi)}\left|\log\frac{v(x)}{w(x)}\right|=d(v,w)\qquad\text{for every }v,w\in\Kcal.
    \end{align*}

    The map \eqref{MapLogSin} is onto. Indeed, for $u\in C_b((0,\pi))$ the function $v:=e^u\sin$ is continuous on $(0,\pi)$ and satisfies
    \begin{align}\label{IneqComparableToSine}
        e^{-\norm{u}_{C_b((0,\pi))}}\sin x\leq v(x)\leq e^{\norm{u}_{C_b((0,\pi))}}\sin x.
    \end{align}
    From \eqref{IneqComparableToSine}, $v(x)\to0$ as $x\to0^+$ and as $x\to\pi^-$, hence $v$ extends continuously to $[0,\pi]$ and $v\in\Kcal$. Finally, since $C_b((0,\pi))$ is complete, so is $\Kcal$.
\end{proof}

\begin{prop}
\label{PropGeneralUniqueness}
    Let $\al>1$ and let $\OpSGen:C([0,\pi])\to C([0,\pi])$ be a linear operator such that $\OpSGen v\geq0$ whenever $v\geq0$. Assume that there exist constants $0<c\leq C$ such that
    \begin{align}\label{HypTwoSided}
        c\,\sin x\leq\OpSGen\big(\sin^{\frac{1}{\al}}\big)(x)\leq C\sin x
        \qquad\text{for every }x\in[0,\pi].
    \end{align}
    Then, the operator $\OpNolinearGen v:=\OpSGen\big(v^{\frac{1}{\al}}\big)$ has exactly one fixed point in $\Kcal$.
\end{prop}

\begin{proof}
    Let $v,w\in C([0,\pi])$ with $w-v\geq0$. Since $t\mapsto t^{\frac{1}{\al}}$ is increasing, $w^{\frac{1}{\al}}-v^{\frac{1}{\al}}\geq0$ too. As $\OpSGen$ is linear and maps nonnegative functions to nonnegative functions,
    \begin{align}\label{IneqMonotone}
        \OpNolinearGen w-\OpNolinearGen v=\OpSGen\big(w^{\frac{1}{\al}}-v^{\frac{1}{\al}}\big)\geq0.
    \end{align}

    Now, notice that $\OpNolinearGen$ maps $\Kcal$ into $\Kcal$. Indeed, if $\xi\sin\leq v\leq\Xi\sin$, then \eqref{HypTwoSided}, \eqref{IneqMonotone} and the linearity of $\OpSGen$ give
    \begin{align*}
        \xi^{\frac{1}{\al}}c\,\sin x\leq
        \xi^{\frac{1}{\al}}\OpSGen\big(\sin^{\frac{1}{\al}}\big)(x)=
        \OpNolinearGen(\xi\sin)(x)\leq
        \OpNolinearGen v(x),
    \end{align*}
    and $\OpNolinearGen v(x)\leq\Xi^{\frac{1}{\al}}C\sin x$ follows similarly.

    Finally, let $v,w\in\Kcal$, so that $e^{-d(v,w)}w\leq v\leq e^{d(v,w)}w$ on $[0,\pi]$. Applying $\OpNolinearGen$ and using \eqref{IneqMonotone},
    \begin{align*}
        \OpNolinearGen v\leq\OpNolinearGen\big(e^{d(v,w)}w\big)=e^{\frac{d(v,w)}{\al}}\OpNolinearGen w,
        \qquad
        \OpNolinearGen w\leq e^{\frac{d(v,w)}{\al}}\OpNolinearGen v,
    \end{align*}
    which is exactly $d(\OpNolinearGen v,\OpNolinearGen w)\leq\tfrac{1}{\al}d(v,w)$. As $\al>1$, the map $\OpNolinearGen$ is a contraction of the complete metric space of Lemma~\ref{LemmaComparableComplete} into itself, and the Banach fixed point theorem \cite{Banach1922} gives the conclusion.
\end{proof}

\subsection{Regularity of the solution}\label{regularity}

This section collects the regularity results used below. Lemma~\ref{LemmaRegularityHolder} deals with a solution $v_\al$ of \eqref{EqnIntegral} on $[0,\pi]$, and Corollary~\ref{CorRegularitySolution} transfers the result to the corresponding solution $f_\al$ of \eqref{EqnStatFolded} on $\T$, which is the form used in Theorems~\ref{TheoremAlBig} and \ref{TheoremAlSmall} below.
\begin{lemma}\label{LemmaRegularityHolder}
    Let $\al>0$ and let $v_\al\in C([0,\pi])$ be a nonnegative solution of \eqref{EqnIntegral} with $n\geq 3$, and denote $g_\al:=-v_\al^{1/\al}$.
    Then
    \begin{align*}
        v_\al\in C^{1,\ga}([0,\pi])
        \quad\text{and}\quad
        g_\al\in C^{\ga}([0,\pi])
        \qquad \text{for every }
        \ga<\min\!\big\{1,\tfrac{1}{\al}\big\}.
    \end{align*}
\end{lemma}
\begin{proof}
    The proof follows the argument in \cite[Proposition~3]{OkaSaWu}, which relies on the boundedness of the Hilbert transform $H$ on $L^p([0,\pi]),\ 1<p<\infty$, and on $C^\ga([0,\pi]),\ 0<\ga<1$, so we will not repeat the full argument here.

    Nevertheless, we work with a different operator, which differs from $\dx^{-1}\!H$ by a term that decays sufficiently fast:
    \begin{align*}
        \OpS=\dx^{-1}\!H(I + \OpSmooth).
    \end{align*}
    Here $\OpSmooth$ is the smoothing Fourier multiplier defined by
    \begin{align*}
        (\ha{\OpSmooth})_k:=\frac{3}{n^2k^2-4},\ k\neq 0.
    \end{align*}

    Since $(\ha{\OpSmooth})_k=O(k^{-2})$, the operator $\OpSmooth$ is bounded on $W^{s,p}([0,\pi]),\ 1<p<\infty$, and on $C^\ga([0,\pi]),\ 0<\ga<1$. Hence the same holds for $\dx\OpS=H(I+\OpSmooth)$, and the bootstrap applies without any change.
\end{proof}

\begin{corollary}\label{CorRegularitySolution}
    Let $\al>\tfrac{1}{2}$ with $\al\neq 1$, let $n\geq 3$, and let $v_\al\in C([0,\pi])$ be a solution of \eqref{EqnIntegral} that is positive on $(0,\pi)$ and vanishes at $0$ and $\pi$. Let $g_\al$ be as in Lemma~\ref{LemmaRegularityHolder} and let $\ti g_\al$ be its odd extension to $\T$, and set $f_\al(x):=\ti g_\al(nx)$. Then $f_\al$ is odd and $\frac{2\pi}{n}$-periodic, it solves \eqref{EqnStatFolded}, and $(\ha{f_\al})_n\neq 0$. Moreover, $f_\al$ belongs to $C^{\frac{1}{\al}}(\T)$ if $\al>1$, and to $C^{1,\frac{1}{\al}-1}(\T)$ if $\al<1$.
\end{corollary}
\begin{proof}
    By Lemma~\ref{LemmaRegularityHolder}, $v_\al\in C^{1,\ga}([0,\pi])$ and $g_\al\in C^{\ga}([0,\pi])$ for every $\ga<\min\!\big\{1,\tfrac{1}{\al}\big\}$. When $\al<1$ we may in addition take $\ga\geq\tfrac{1}{\al}-1$, since $\al>\tfrac{1}{2}$ gives $\tfrac{1}{\al}-1<1$.

    Since $\OpS g_\al=v_\al$ by \eqref{EqnIntegral}, these are the hypotheses of Lemma~\ref{FromFixedPointToEQN}, which therefore gives that $\ti g_\al$ belongs to $C^{\frac{1}{\al}}(\T)$ if $\al>1$ and to $C^{1,\frac{1}{\al}-1}(\T)$ if $\al<1$, and that it solves \eqref{EqnStatFoldedg}.

    The function $f_\al(x):=\ti g_\al(nx)$ is odd and $\frac{2\pi}{n}$-periodic, and it belongs to the same space as $\ti g_\al$, since rescaling preserves the H\"older exponent. By Lemma~\ref{Lemma_FourierUnfoldFold}, $f_\al$ solves \eqref{EqnStatFolded} and $(\ha{f_\al})_n$ is a nonzero multiple of the first sine coefficient of $g_\al$. Since $v_\al=\OpS g_\al$ and $(\ha{\OpS})_1\neq 0$, that coefficient is a nonzero multiple of $\int_0^\pi v_\al(x)\sin x\,dx$.
\end{proof}

\subsection{Verifying Schauder's hypotheses for $-\OpS$}\label{checking}
In this section we prove that $-\OpS$, with $\OpS$ defined in \eqref{EqnStatUnFolded}, satisfies the assumptions of Proposition~\ref{PropGeneralExistence} with $\OpSGen=-\OpS$ when $\al\in\big(\tfrac{1}{2},1\big)$, where, as above, $\OpS$ acts on a function on $[0,\pi]$ by applying it to the odd extension and restricting back to $[0,\pi]$.

A direct computation gives
\begin{align}
    \int_0^\pi
    -\OpS\big(v^\frac{1}{\al}\big)(x)\sin x\,dx&=
    \frac{2}{\pi}\int_0^\pi\bigg(\int_0^\pi
    \sum_{k\geq 1}
    \frac{1}{k}\frac{n^2k^2-1}{n^2k^2-4}
    v^\frac{1}{\al}(y)
    \sin(kx)\sin(ky)\,dy\bigg)\sin x\,dx\nonumber\\
    &=
    \frac{n^2-1}{n^2-4}\int_0^\pi
    v^\frac{1}{\al}(y)
    \sin y\,dy.\label{a}
\end{align}
Identity \eqref{a} corresponds to property~\ref{enume_first} in Proposition~\ref{PropGeneralExistence} with $\rho=\tfrac{n^2-1}{n^2-4}$.

Moreover, since $(\ha{\OpS})_k$ behaves like $-\frac{1}{|k|}$ as $k$ goes to infinity, $\OpS$ is compact on $C([0,\pi])$.

To verify hypothesis~\ref{enume_second} in Proposition~\ref{PropGeneralExistence}, we next show,
following the strategy of \cite{OkaSaWu}, that $-\OpS\big(v^\frac{1}{\al}\big)$ admits an integral representation
on $[0,\pi]$ against a certain singular kernel.
\begin{lemma}\label{Lemma_Kernel_Def}
    Let $v\in C([0,\pi])$ be nonnegative and $n\geq 3$. Then, for any $x\in[0,\pi]$ we have that
    \begin{align*}
        -\OpS\big(v^{\frac{1}{\al}}\big)(x)
        =
        \int_0^\pi
        v(y)^\frac{1}{\al}\,
        \kernel(x,y)\,dy,
    \end{align*}
    where
    \begin{align*}
        \kernel(x,y)
        :=
        \frac{1}{2\pi}\log
        \!\left(\frac{1-\cos(x+y)}{1-\cos(x-y)}\right)
        +
        \frac{6}{\pi}\sum_{k\geq 1}
        \frac{1}{n^2k^3-4k}
        \sin(kx)\sin(ky).
    \end{align*}
\end{lemma}
\begin{proof}
    We write $\OpS=\dx^{-1}\!H(I+\OpSmooth)$, where the multiplier of $\OpSmooth$ is given by
    \begin{align*}
        (\ha{\OpSmooth})_k:=
        \frac{3}{n^2k^2-4}.
    \end{align*}
    We can write $\dx^{-1}\!H\big(v^{\frac{1}{\al}}\big)$ as
    \begin{align*}
        \dx^{-1}\!H\big(v^\frac{1}{\al}\big)(x) = \frac{1}{2\pi}\intpi v(y)^\frac{1}{\al}\log(1-\cos(x-y)) \,dy.
    \end{align*}
    Assuming now that $v$ is odd, we obtain
    \begin{align*}
        \dx^{-1}\!H\OpSmooth \big(v^\frac{1}{\al}\big)(x) &=
        -\frac{6}{\pi}
        \sum_{k\geq 1}
        \int_0^\pi
        v(y)^\frac{1}{\al}
        \frac{1}{n^2k^3-4k}
        \sin(kx)\sin(ky)\,dy,
        \\
        \dx^{-1}\!H\big(v^\frac{1}{\al}\big)(x) &=
        \frac{1}{2\pi}\int_0^\pi v(y)^{\frac{1}{\al}}
        \log\!
        \left(\frac{1-\cos(x-y)}{1-\cos(x+y)}\right) \,dy.
    \end{align*}
    Thus,
    \begin{align*}
        -\OpS \big(v^\frac{1}{\al}\big)(x)
        &=
        -\dx^{-1}\!H\big(v^\frac{1}{\al}\big)(x)
        -\dx^{-1}\!H\OpSmooth \big(v^\frac{1}{\al}\big)(x)
        \\ &=
        \frac{1}{2\pi}\int_0^\pi v(y)^\frac{1}{\al}
        \log\!
        \left(\frac{1-\cos(x+y)}{1-\cos(x-y)}\right)
        \, dy
        \\
        &\qquad+
        \frac{6}{\pi}\sum_{k\geq 1}\int_0^\pi
        v(y)^\frac{1}{\al}
        \frac{1}{n^2k^3-4k}
        \sin(kx)\sin(ky)
        \, dy\\
        &=
        \int_0^\pi
        v(y)^\frac{1}{\al}\,
        \kernel (x,y)\,dy.
    \end{align*}
\end{proof}
We will use the following three lemmas to prove that $-\OpS\big(v^{\frac{1}{\al}}\big) \geq 0$ for every nonnegative $v\in C([0,\pi])$.
\begin{lemma}\label{LemmaSineSeriesPositive}
    For every even $s\geq 2$, let $f_s$ be the function given by
    \begin{align*}
        f_s(x):=\sum_{k\geq 1}\frac{\sin(kx)}{k^s},
    \end{align*}
    so that $f_2$ is the Clausen function $\clau$. Then, $f_s(x)\geq 0$ for every $x\in[0,\pi]$.
\end{lemma}
\begin{proof}
    Each of these series converges absolutely, so $f_s$ is continuous on $[0,\pi]$, and $f_s(0)=f_s(\pi)=0$ because $\sin(k\pi)=0$ for every $k\geq 1$.

    For $s=2$ we use the expansion
    \begin{align}
        \log(2-2\cos x)&=
        \log(1-e^{ix})+\log(1-e^{-ix})\nonumber\\&=
        -\sum_{k\geq 1} \frac{e^{ikx}}{k}
        -\sum_{k\geq 1} \frac{e^{-ikx}}{k}=
        -2\sum_{k\geq 1} \frac{1}{k}\cos(kx),\label{equ_expan_log}
    \end{align}
    whose right hand side converges uniformly on compact subsets of $(0,2\pi)$. Differentiating the series of $\clau$ term by term is therefore legitimate there and gives
    \begin{align}\label{clauder}
        \dx\clau(x)=-\frac{1}{2}\log(2-2\cos x),
    \end{align}
    so that $\dx\clau(x)>0$ if and only if $2-2\cos x<1$, that is, if and only if $x<\tfrac{\pi}{3}$. Hence $\clau$ increases on $\left(0,\tfrac{\pi}{3}\right)$ and decreases on $\left(\tfrac{\pi}{3},\pi\right)$. Since it vanishes at both endpoints, it is nonnegative on $[0,\pi]$.

    Assume now that $f_s\geq 0$ on $[0,\pi]$ for some even $s\geq2$. Differentiating the series of $f_{s+2}$ twice term by term produces the series of $-f_s$, which converges uniformly because $s\geq2$. Hence $f_{s+2}\in C^2([0,\pi])$ and
    \begin{align*}
        \dx^2 f_{s+2}=-f_s\leq 0\quad\text{on }[0,\pi],
    \end{align*}
    so, since $f_{s+2}$ also vanishes at both endpoints of $[0,\pi]$, it is nonnegative there.
\end{proof}
\begin{lemma}\label{LemmaPhiIncreasing}
    Let $n\geq3$ and let $\phi_n$ be the function given by
    \begin{align}\label{DefPhiN}
        \phi_n(x):=\log(1-\cos x)+\psi_n(x),
        \qquad
        \psi_n(x):=-6\sum_{k\geq 1}\frac{\cos(kx)}{n^2k^3-4k}.
    \end{align}
    Then, $\dx\psi_n\geq0$ on $[0,\pi]$ and $\dx\phi_n>0$ on $(0,\pi)$. In particular, both $\psi_n$ and $\phi_n$ are nondecreasing on $[0,\pi]$, and $\phi_n$ is strictly increasing there.
\end{lemma}
\begin{proof}
    The series defining $\psi_n$ and the one obtained by differentiating it term by term both converge uniformly, so
    \begin{align}\label{DerPsiN}
        \dx\psi_n(x)=6\sum_{k\geq 1}\frac{\sin(kx)}{n^2k^2-4}.
    \end{align}
    Set $a:=\tfrac{2}{n}$, so that $0<a\leq\tfrac{2}{3}$ for $n\geq3$. Since $a^2k^{-2}\leq\tfrac{4}{9}<1$ for every $k\geq1$, expanding a geometric series gives
    \begin{align}\label{EqnGeomExpansion}
        \frac{1}{n^2k^2-4}=
        \frac{1}{n^2}\,\frac{1}{k^2-a^2}=
        \frac{1}{n^2}\sum_{j\geq 0}\frac{a^{2j}}{k^{2j+2}}.
    \end{align}
    Since every summand in \eqref{EqnGeomExpansion} is positive, and $\sum_{k\geq1}(k^2-a^2)^{-1}<\infty$, the resulting double series converges absolutely, so it can be rearranged freely. Substituting \eqref{EqnGeomExpansion} into \eqref{DerPsiN} and rearranging,
    \begin{align*}
        \dx\psi_n(x)=\frac{6}{n^2}\sum_{j\geq 0}a^{2j}f_{2j+2}(x),
    \end{align*}
    Here $f_{2j+2}$ denotes the function from Lemma~\ref{LemmaSineSeriesPositive}, which is nonnegative on $[0,\pi]$ for every $j\geq0$ by that lemma. Since $a\geq0$, this shows $\dx\psi_n\geq0$ on $[0,\pi]$. Finally, $\dx\log(1-\cos x)=\cot\big(\frac x2\big)>0$ for $x\in(0,\pi)$, so $\dx\phi_n>0$ there.
\end{proof}
\begin{lemma}\label{LemmaEvenMonotone}
    Let $g:\RR\to\RR\cup\{-\infty\}$ be even, $2\pi$-periodic and nondecreasing on $[0,\pi]$. Then,
    \begin{align*}
        g(x+y)\geq g(x-y)
        \qquad\text{for every }(x,y)\in[0,\pi]^2.
    \end{align*}
    If moreover $g$ is strictly increasing on $[0,\pi]$, the inequality is strict whenever $x,y\in(0,\pi)$.
\end{lemma}
\begin{proof}
    We work with the variables $u:=x+y$ and $z:=x-y$. As $(x,y)$ ranges over $[0,\pi]^2$, the pair
    $(u,z)$ ranges over the square
    \begin{align*}
        R_{\text{square}}:=\{(u,z):\ 0\leq u+z\leq 2\pi,\ 0\leq u-z\leq2\pi\},
    \end{align*}
    whose vertices are $(0,0)$, $(\pi,\pi)$, $(2\pi,0)$ and $(\pi,-\pi)$. In these variables the
    statement reads $g(u)\geq g(z)$ for every $(u,z)\in R_{\text{square}}$.

    We first use the symmetries of $g$ to reduce the problem to the triangle
    \begin{align*}
        R_{\text{triangle}}:=\{(u,z):\ 0\leq z\leq u\leq \pi\}.
    \end{align*}
    Since $g(z)=g(-z)$ and $g(u)=g(2\pi-u)$, we can assume that $(u,z)\in R_{\text{triangle}}$
    rather than merely $R_{\text{square}}$, see Figure~\ref{fig:RsquareRtriangle}. It is therefore
    enough to prove the inequality on $R_{\text{triangle}}$.

    \begin{figure}[H]
        \centering
    \begin{tikzpicture}[scale=1]
        \pgfmathsetmacro{\xmin}{-1}
        \pgfmathsetmacro{\xmax}{2*pi+1}
        \pgfmathsetmacro{\ymin}{-1}
        \pgfmathsetmacro{\ymax}{pi}

        \tkzDefPoint(0,0){O}
        \tkzDefPoint(pi,0){A}
        \tkzDefPoint(pi,pi){B}
        \tkzDefPoint(pi,-pi){C}
        \tkzDefPoint(2*pi,0){D}

        \draw[->,thin] (-1,0) -- (1,0) node[below] {$u$};
        \draw[->,thin] (0,-1) -- (0,1) node[below right] {$z$};

        \draw[-,dash dot] (\xmin-1,0) -- (\xmax+1,0);
        \draw[-,dash dot] (pi,\ymin-pi) -- (pi,\ymax+1);

        \tkzDrawPolygon[thick](O,B,D,C)

        \tkzFillPolygon[fill=black!10,opacity=0.6](O,A,B)
        \tkzDrawPolygon[thick](O,A,B)

        \tkzDrawPoints[fill=black,size=6](O,A,B,C,D)

        \tkzLabelPoint[below left](O){$(0,0)$}
        \tkzLabelPoint[below right](A){$\!(\pi,0)$}
        \tkzLabelPoint[above right](B){$\!(\pi,\pi)$}
        \tkzLabelPoint[below right](C){$(\pi,-\pi)$}
        \tkzLabelPoint[below right](D){$\!(2\pi,0)$}

        \node at (2.13,0.9) {$R_{\text{triangle}}$};
        \node at (5.5,1.7) {$R_{\text{square}}$};

        \node at (7.9,0.35) {$g(-z)=g(z)$};
        \node at (1.35,-3.7) {$g(2\pi-u)=g(u)$};
    \end{tikzpicture}
    \caption{Representation of $R_{\text{square}}$ and $R_{\text{triangle}}$.}
    \label{fig:RsquareRtriangle}
    \end{figure}
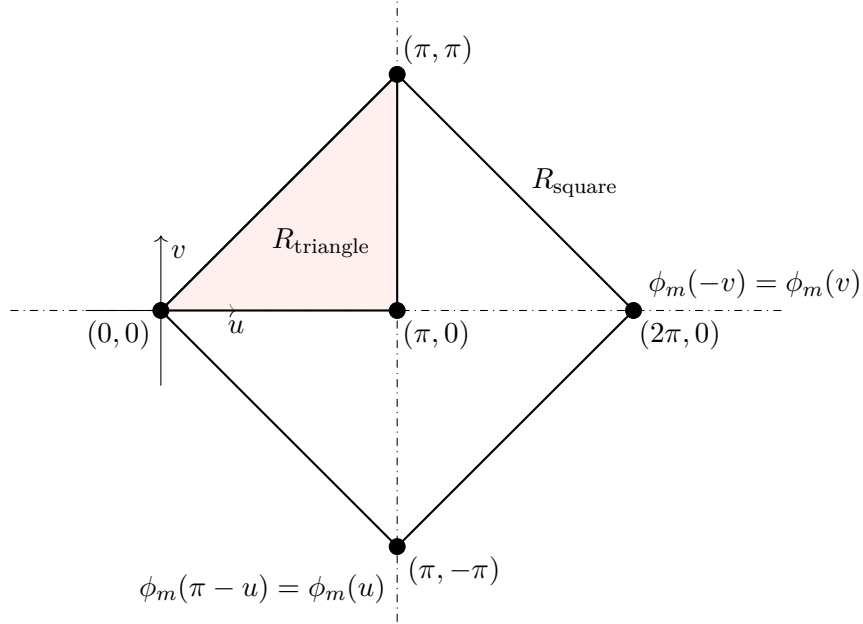

    Indeed, if $0\leq z\leq u\leq\pi$, then $g(u)\geq g(z)$ since $g$ is nondecreasing there.

    Assume now that $g$ is strictly increasing on $[0,\pi]$. The inequality on $R_{\text{triangle}}$
    is then strict as soon as $z<u$, so equality is achieved only when $z=u$. Since
    $(x,y)\mapsto(u,z)$ maps the boundary of $[0,\pi]^2$ to the boundary of $R_{\text{square}}$, the
    inequality of the statement is strict for every $(x,y)\in(0,\pi)^2$.
\end{proof}
The following lemma shows that the kernel $\kernel$ is nonnegative, which gives the first inequality of hypothesis~\ref{enume_second} in Proposition~\ref{PropGeneralExistence}.
\begin{lemma}\label{LemmaKernelPositive}
    Let $n\geq 3$ and let $v\in C([0,\pi])$ be nonnegative and not identically zero. Then,
    \begin{align*}
        \kernel(x,y)\geq 0\quad\text{for a.e. }(x,y)\in [0,\pi]^2,
    \end{align*}
   and $-\OpS\big(v^{\frac{1}{\al}}\big)\in C([0,\pi])$ is nonnegative.
   Furthermore,
   \begin{align*}
       -\OpS\big(v^{\frac{1}{\al}}\big)(x)>0\quad\text{for every }x\in(0,\pi),\quad -\OpS\big(v^{\frac{1}{\al}}\big)(0)=-\OpS\big(v^{\frac{1}{\al}}\big)(\pi)=0.
   \end{align*}
\end{lemma}
\begin{proof}
    Since the kernel $\kernel$ is integrable, $-\OpS\big(v^{\frac{1}{\al}}\big)$ is continuous on $[0,\pi]$ for all $v\in C([0,\pi])$.
    In addition, $\kernel$ satisfies
    \begin{align*}
        \kernel(0,y)=\kernel(\pi,y)=0\quad\text{for all }y\in[0,\pi].
    \end{align*}
    Hence, $-\OpS\big(v^{\frac{1}{\al}}\big)(0)=-\OpS\big(v^{\frac{1}{\al}}\big)(\pi)=0$.

    The nontrivial part of the statement is the nonnegativity. Recall from Lemma~\ref{Lemma_Kernel_Def} that
    \begin{align*}
        -\OpS\big(v^{\frac{1}{\al}}\big)(x)=
        \int_0^\pi
        v(y)^{\frac{1}{\al}}\,
        \kernel(x,y)\,dy.
    \end{align*}
    The kernel $\kernel(x,y)$ is not defined on the diagonal $x=y$, but this set has measure zero and the singularity is integrable. Hence, it is enough to prove that $\kernel(x,y)\geq 0$ for a.e. $(x,y)\in[0,\pi]^2$.
    We work with the variables $u:=x+y$ and $z:=x-y$, which lie in
    \begin{align*}
        R_{\text{square}}:=\{(u,z):\ 0<u+z<2\pi,\ 0<u-z<2\pi\}.
    \end{align*}
    Then, we obtain
    \begin{align*}
        \ti{K}_n(u,z) &=
        \kernel\!\left(\frac{u+z}{2},\frac{u-z}{2}\right)\\ &=
        \frac{1}{2\pi}\log\!\left(\frac{1-\cos u}{1-\cos z}\right) +
        \frac{3}{\pi}\sum_{k\geq 1}
        \frac{1}{n^2k^3-4k}\left(\cos(kz) - \cos(ku)\right)\\ &=
        \frac{1}{2\pi}\phi_n(u) - \frac{1}{2\pi}\phi_n(z).
    \end{align*}
    Here, $\phi_n$ is the function introduced in Lemma~\ref{LemmaPhiIncreasing}. Since $\phi_n$ is even and $2\pi$-periodic, and $\ti{K}_n(u,z)=\frac{1}{2\pi}\phi_n(u)-\frac{1}{2\pi}\phi_n(z)$ for every $(u,z)\in R_{\text{square}}$, we obtain the symmetries
    \begin{align}\label{symWiKernel}
        \ti{K}_n(u,z)=
        \ti{K}_n(u,-z)=
        \ti{K}_n(2\pi-u,z)=
        \ti{K}_n(2\pi-u,-z).
    \end{align}
    These symmetries will be used later, in the proof of Lemma~\ref{LemmaKernelBoundSine}.

    By Lemma~\ref{LemmaPhiIncreasing}, $\phi_n$ is strictly increasing on $[0,\pi]$, so Lemma~\ref{LemmaEvenMonotone} applied to $g=\phi_n$ gives
    \begin{align*}
        \kernel(x,y)=\frac{1}{2\pi}
        \big(\phi_n(x+y)-\phi_n(x-y)\big)\geq 0
        \qquad\text{for a.e. }(x,y)\in[0,\pi]^2,
    \end{align*}
    with strict inequality whenever $x,y\in(0,\pi)$.

    Finally, fix $x\in(0,\pi)$. Then $\kernel(x,y)>0$ for every $y\in(0,\pi)$, and $v$ is nonnegative and does not vanish identically, so
    \begin{align*}
        -\OpS\big(v^{\frac{1}{\al}}\big)(x)=
        \int_0^\pi
        v(y)^{\frac{1}{\al}}\,
        \kernel(x,y)\,dy>0.
    \end{align*}
\end{proof}
The following lemma proves the second inequality of hypothesis~\ref{enume_second} in Proposition~\ref{PropGeneralExistence}.
\begin{lemma}\label{LemmaBoundAlpha}
    Let $\al\in\big(\tfrac{1}{2},1\big)$. Then, there exists $M>0$ such that $\ConvSet\neq\emptyset$ and, for every $v\in\ConvSet$, the following estimate holds,
    \begin{align*}
        \norm{-\OpS\big(v^\frac{1}{\al}\big)}_{C([0,\pi])}\leq
        M\int_0^\pi -\OpS \big(v^{\frac{1}{\al}}\big)(y)\,\sin y\,dy.
    \end{align*}
\end{lemma}

\begin{rem}
    We cannot choose $M$ arbitrarily small. Indeed, if $M<\tfrac{1}{2}$ and $v\in\ConvSet$, then
    \begin{align*}
    1 = \int_{0}^\pi v(x)\,\sin x\,dx
    \le M \int_{0}^\pi \sin x\,dx
    = 2M < 1,
    \end{align*}
    a contradiction. Hence, $\ConvSet=\emptyset$ whenever $M<\tfrac{1}{2}$.
\end{rem}
The proof of Lemma~\ref{LemmaBoundAlpha} relies on the following four auxiliary lemmas, whose proofs are deferred to the appendix.
\begin{lemma}\label{LemmaKernelBoundSine}
    For every $(x,y)\in[0,\pi]^2$, it holds that
    \begin{align*}
        \kernel(x,y)\lesssim \sin y\,\cot\!\left(\frac{|x-y|}{2}\right)+\sin y.
    \end{align*}
\end{lemma}
\begin{lemma}\label{Lemma_when_isnegative}
    Let $A_1,A_2>0,\ \al\in(0,1)$ and let $\beta:(0,\infty)\to\RR$ be the function given by
    \begin{align*}
        \beta(M):=
        A_1 M^{\frac{1}{\al}} + A_2 - M.
    \end{align*}
    Then, $\beta$ attains a negative minimum at $M_0:=(A_1\al^{-1})^{\frac{1}{1-\al^{-1}}}$ if and only if
    \begin{align*}
        A_1^\al A_2^{1-\al} < (1-\al)^{1-\al} \al^\al.
    \end{align*}
\end{lemma}
\begin{lemma}
\label{LemmaLimitAlSmallIsZero}
Let $\al\in\big(\tfrac{1}{2},1\big)$ and let $\Ical:\RR^+\to\RR$ be given by
\begin{align*}
    \Ical(\de):=\sup_{x\in[0,\pi]}
    \int_0^\pi \kernel(x,y)\,
    \chi_{|x-y|<\de}(x,y)\,dy.
\end{align*}
Then, we have that
\begin{align*}
\lim_{\de\to0^+}\
\Ical(\de)^\al
\left(1+\cot(\de/2)\right)^{1-\al}=0.
\end{align*}
\end{lemma}
\begin{lemma}
\label{LemmaEpDeltaGood}
Let $\Ical:\RR^+\to\RR$ be the function defined in Lemma~\ref{LemmaLimitAlSmallIsZero}. Then, for every $\de\in\big(0,\tfrac{\pi}{3}\big)$, we have that
\begin{align*}
    \Ical(\de) \lesssim \de(1-\log\de).
\end{align*}
\end{lemma}
\begin{proofsth}{Lemma~\ref{LemmaBoundAlpha}}
    By Lemma~\ref{LemmaKernelPositive}, $-\OpS\big(v^{\frac{1}{\al}}\big)$ is continuous, so its $C([0,\pi])$ norm is well defined. Using identity \eqref{a}, the statement of this lemma is equivalent to the estimate
    \begin{align}\label{IneqWeProveSineBound}
        \norm{-\OpS\big(v^\frac{1}{\al}\big)}_{C([0,\pi])}\leq
        M\frac{n^2-1}{n^2-4}\int_0^\pi v(y)^{\frac{1}{\al}}\,\sin y\,dy.
    \end{align}

We split the integration domain into two parts
\begin{align}
    -\OpS\big(v^{\frac{1}{\al}}\big)(x)
    &=
    \int_0^\pi
    v(y)^\frac{1}{\al}\,
    \kernel(x,y)\,
    dy
    \nonumber\\
    &=\int_0^\pi
    v(y)^\frac{1}{\al}\,
    \kernel(x,y)\,
    \chi_{|x-y|<\de}(x,y)\,
    dy+
    \int_0^\pi
    v(y)^\frac{1}{\al}\,
    \kernel(x,y)\,
    \chi_{|x-y|>\de}(x,y)\,
    dy. \label{eqnSplit}
\end{align}
Here, $\de>0$ will be chosen later.

Recall $\Ical(\de)$ defined in Lemma~\ref{LemmaLimitAlSmallIsZero}. Now, for any $v \in \ConvSet$, we have that
\begin{align}\label{IneqIcalEstimate}
    \int_0^\pi
    v(y)^\frac{1}{\al}\,
    \kernel(x,y)\,
    \chi_{|x-y|<\de}(x,y)\,dy\leq
    M^{\frac{1}{\al}}
    \int_0^\pi \kernel(x,y)\,
    \chi_{|x-y|<\de}(x,y)\,dy\leq
    M^{\frac{1}{\al}}\,\Ical(\de).
\end{align}

Denote now by $\mathcal{B}(\de,x,v)$ the quantity given by
\begin{align*}
    \mathcal{B}(\de,x,v)&:=
    \int_0^\pi
    v(y)^\frac{1}{\al}\,
    \kernel(x,y)\,
    \chi_{|x-y|>\de}(x,y)\,dy.
\end{align*}
By Lemma~\ref{LemmaKernelBoundSine}, we have
\begin{align}
    \mathcal{B}(\de,x,v)&\lesssim
    \left(\cot(\de/2)+1\right)
    \int_0^\pi
    v(y)^{\frac{1}{\al}}\,
    \sin y\,dy.\label{BBound}
\end{align}

For $v \in \ConvSet$, we recall that inequality \eqref{IneqDevelLowerBoundAlSma} can be rewritten as
    \begin{align}\label{IneqLowerBoundAlSma}
        2^{1-\frac{1}{\al}}\leq \int_0^\pi v(x)^{\frac{1}{\al}}\sin x\,dx.
    \end{align}
    Substituting \eqref{IneqLowerBoundAlSma} into \eqref{IneqIcalEstimate} yields
    \begin{align}\label{IneqAlSmaIcal}
        \int_0^\pi
        v(y)^{\frac{1}{\al}}\,
        \kernel(x,y)\,
        \chi_{|x-y|<\de}(x,y)\,dy\leq
        2^{\frac{1}{\al}-1}\,
        M^{\frac{1}{\al}}\,\Ical(\de)\,\int_0^\pi v(x)^{\frac{1}{\al}}\,\sin x\,dx.
    \end{align}
    Consequently, substituting estimates \eqref{BBound} and \eqref{IneqAlSmaIcal} into \eqref{eqnSplit} yields
    \begin{align}\label{eqnEstimateIcalCot}
        -\OpS\big(v^\frac{1}{\al}\big)(x)\lesssim
        \left(M^{\frac{1}{\al}}\Ical(\de)
        +\cot(\de/2)+1\right)
        \int_0^\pi v(y)^{\frac{1}{\al}}\,\sin y\,dy.
    \end{align}
    Thus, it is enough to prove that, for any $\ep>0$, there exist $\de$ and $M$ such that
    \begin{align}\label{dmenor}
        \Ical(\de)M^{\frac{1}{\al}}+1+\cot(\de/2)\leq \ep M.
    \end{align}

    Notice that \eqref{dmenor} is equivalent to showing that, for any $\ep>0$, there exist $\de$ and $M$ such that
    \begin{align*}
        \underbrace{\frac{\Ical(\de)}{\ep}}_{A_1}M^\frac{1}{\al}+
        \underbrace{\frac{1+\cot\!\left(\de/2\right)}{\ep}}_{A_2}
        -M\leq 0.
    \end{align*}
    By Lemma~\ref{Lemma_when_isnegative}, this is equivalent to proving that, for any $\ep>0$, there exists $\de>0$ such that
    \begin{align*}
        \left(\frac{\Ical(\de)}{\ep}\right)^\al
        \left(\frac{1+\cot(\de/2)}{\ep}\right)^{1-\al}
        =
        A_1^{\al} A_2^{1-\al}
        < (1-\al)^{1-\al} \al^\al.
    \end{align*}
    This is possible if we can make the product
    \begin{align*}
        \Ical(\de)^\al
        \left(1+\cot\!\left(\de/2\right)\right)^{1-\al}
    \end{align*}
        arbitrarily small by taking $\de$ small enough, which is true since Lemma~\ref{LemmaLimitAlSmallIsZero} proves that, for $\al\in\big(\tfrac{1}{2},1\big)$,
    \begin{align*}
        \lim_{\de\to0^+}\,
        \Ical(\de)^\al
        \left(1+\cot\!\left(\de/2\right)\right)^{1-\al}=0.
    \end{align*}

    By Lemma~\ref{Lemma_when_isnegative}, we take $M$ to be
    \begin{align*}
        M_0:=\left(\frac{\Ical(\de)}{\ep}\al^{-1}\right)^{\frac{1}{1-\al^{-1}}}.
    \end{align*}
    Notice that taking $\de$ smaller gives a larger $M_0$, so that $\ConvSet$ is never empty.
\end{proofsth}

We have proved that $-\OpS$ satisfies every assumption in Proposition~\ref{PropGeneralExistence} with $\OpSGen=-\OpS$ and $\al\in\big(\tfrac{1}{2},1\big)$.

\subsection{Verifying Banach's hypothesis for $-\OpS$}\label{twosided}
We now check hypothesis \eqref{HypTwoSided} of Proposition~\ref{PropGeneralUniqueness} for $-\OpS$. We start with the logarithmic term of $\kernel$.
\begin{lemma}\label{LemmaKernelTwoSided}
    Let $\al>0$. Then, there exist constants $0<c\leq C$, depending only on $\al$, such that
    \begin{align*}
        c\,\sin x\leq
        \frac{1}{2\pi}
        \int_0^\pi
        \sin^{\frac{1}{\al}}y\,
        \log\!\left(\frac{1-\cos(x+y)}{1-\cos(x-y)}\right)dy
        \leq C\,\sin x
        \quad\text{for every }x\in[0,\pi].
    \end{align*}
\end{lemma}
\begin{proof}
    For $x\in(0,\pi)$, let $h(x)$ be the quotient of the middle term above by $\sin x$. We claim that $h$ extends to a positive continuous function on $[0,\pi]$, which proves the lemma with $c:=\min h$ and $C:=\max h$.

    Write $s:=\tan\big(\frac{x}{2}\big)$ and $u:=\tan\big(\frac{y}{2}\big)$, so that
    \begin{align}\label{IdentityTangentKernel}
        \log\!\left(\frac{1-\cos(x+y)}{1-\cos(x-y)}\right)
        =
        2\log\left|\frac{s+u}{s-u}\right|.
    \end{align}
    Therefore, changing variables we get
    \begin{align}
        h(x)&=
        \frac{2^{\frac{1}{\al}}}{\pi}\,
        \frac{1+s^2}{s}
        \int_0^\infty
        w(u)\,
        \log\left|\frac{s+u}{s-u}\right|du,
        \qquad\text{with}\qquad
        w(u):=\frac{u^{\frac{1}{\al}}}{\left(1+u^2\right)^{\frac{1}{\al}+1}},
        \nonumber\\
        &=
        \frac{2^{\frac{1}{\al}}}{\pi}\big(1+s^2\big)
        \int_0^\infty w(st)\,\log\left|\frac{1+t}{1-t}\right|dt,
        \label{IdentityTwoSidedScaled}
    \end{align}
    Notice that
    \begin{align}\label{weightineq}
        w(u)\leq\min\big\{u^{\frac{1}{\al}},\,u^{-\frac{1}{\al}-2}\big\},
    \end{align}
    which follows from $(1+u^2)^{-1}\leq\min\{1,u^{-2}\}$.

    For $s$ in a compact subset of $(0,\infty)$, from \eqref{weightineq} in the right hand side of \eqref{IdentityTwoSidedScaled} we get
    \begin{align}\label{IneqDominantBound}
        w(st)\,\log\left|\frac{1+t}{1-t}\right|
        \lesssim
        \min\big\{t^{\frac{1}{\al}},\,t^{-\frac{1}{\al}-2}\big\}\log\left|\frac{1+t}{1-t}\right|,
    \end{align}
    with implicit constant not depending on $s$ in the compact subset.

    In order to show that the right hand side of \eqref{IneqDominantBound} is integrable, we have to control its behavior at $0$, $1$ and $\infty$. Notice that
    \begin{align}\label{IneqTwoSidedLogUpper}
        \log\!\left(\frac{1+z}{1-z}\right)=2\sum_{k\geq0}\frac{z^{2k+1}}{2k+1}
        \leq2\sum_{k\geq0}z^{2k+1}
        =\frac{2z}{1-z^2}
        \leq\frac{8}{3}z
        \quad\text{for }0<z\leq\tfrac{1}{2}.
    \end{align}
    Hence, at $0$, recalling \eqref{IneqDominantBound} and taking $z=t$ in \eqref{IneqTwoSidedLogUpper} gives, for $0<t\leq\tfrac{1}{2}$,
    \begin{align}\label{IneqDominantNearZero}
        \min\big\{t^{\frac{1}{\al}},\,t^{-\frac{1}{\al}-2}\big\}\log\left|\frac{1+t}{1-t}\right|
        \leq
        t^{\frac{1}{\al}}\log\!\left(\frac{1+t}{1-t}\right)
        \leq
        \frac{8}{3}\,t^{\frac{1}{\al}+1},
    \end{align}
    while, for $\infty$, taking $z=\tfrac{1}{t}$ in \eqref{IneqTwoSidedLogUpper} gives, for $t\geq2$,
    \begin{align}\label{IneqDominantNearInfty}
        \min\big\{t^{\frac{1}{\al}},\,t^{-\frac{1}{\al}-2}\big\}\log\left|\frac{1+t}{1-t}\right|
        \leq
        t^{-\frac{1}{\al}-2}\log\!\left(\frac{1+t^{-1}}{1-t^{-1}}\right)
        \leq
        \frac{8}{3}\,t^{-\frac{1}{\al}-3}.
    \end{align}
    Both bounds are integrable since $\tfrac{1}{\al}+1>-1$ and $-\tfrac{1}{\al}-3<-1$. At $t=1$ the singularity is only logarithmic, hence also integrable.

    On every compact subset of $(0,\infty)$, the left hand side of \eqref{IneqDominantBound} depends continuously on $s$ for each fixed $t$ and is dominated by an integrable function of $t$. Hence, dominated convergence shows that $h$ is continuous on $(0,\pi)$.

    It remains to compute the limit of $h$ as $s\to0^+$. Since the cosine function is even and $\sin(\pi-y)=\sin y$, the change of variable $y\mapsto\pi-y$ gives
    \begin{align*}
        h(\pi-x)=h(x),
    \end{align*}
    so proving continuity at $x=0$ also gives continuity at $x=\pi$. We split the integral in \eqref{IdentityTwoSidedScaled} at $t=2$. For the part over $(0,2)$, recalling \eqref{weightineq}, we get
    \begin{align}\label{IneqFirstPartVanishes}
        \int_0^2 w(st)\,\log\left|\frac{1+t}{1-t}\right|dt
        \leq
        (2s)^{\frac{1}{\al}}\int_0^2\log\left|\frac{1+t}{1-t}\right|dt
        \lesssim
        s^{\frac{1}{\al}}
        \to0
        \quad\text{as }s\to0^+.
    \end{align}

    For the part over $(2,\infty)$, we go back to the variable $u=st$, since $\tfrac{s}{u}<\tfrac{1}{2}$, from \eqref{IneqTwoSidedLogUpper} we find that
    \begin{align*}
        \frac{1}{s}\int_{2s}^\infty w(u)\,\log\left|\frac{u+s}{u-s}\right|du
        \leq
        \frac{8}{3}\int_0^\infty\frac{w(u)}{u}\,du.
    \end{align*}
    Now notice that this last integral is finite. Indeed, from \eqref{weightineq}, $\frac{w(u)}{u}$ is integrable near $u=0$ since $\tfrac{1}{\al}-1>-1$ and at infinity since $-\tfrac{1}{\al}-3<-1$. As $s\to0^+$, $\frac{1}{s}\,w(u)\,\log\big(\frac{u+s}{u-s}\big)$ converges to $\frac{2w(u)}{u}$ for every fixed $u>0$, and is dominated by the integrable function $\frac{8}{3}\frac{w(u)}{u}$, while $2s<u$. Applying again the dominated convergence theorem gives
    \begin{align}\label{IdentityLimitSecondPart}
        \lim_{s\to0^+}\frac{1}{s}\int_{2s}^\infty w(u)\,\log\left|\frac{u+s}{u-s}\right|du
        =
        \int_0^\infty\frac{2w(u)}{u}\,du.
    \end{align}

    Finally, from \eqref{IdentityTwoSidedScaled}, \eqref{IneqFirstPartVanishes} and \eqref{IdentityLimitSecondPart}, we have
    \begin{align*}
        \lim_{x\to0^+}h(x)
        &=
        \lim_{s\to0^+}\frac{2^{\frac{1}{\al}}}{\pi}(1+s^2)
        \left(\int_0^2 w(st)\,\log\left|\frac{1+t}{1-t}\right|dt+\int_2^\infty w(st)\,\log\left|\frac{1+t}{1-t}\right|dt\right)\\
        &=
        \lim_{s\to0^+}\frac{2^{\frac{1}{\al}}}{\pi}(1+s^2)\cdot\frac{1}{s}\int_{2s}^\infty w(u)\,\log\left|\frac{u+s}{u-s}\right|du\\
        &=
        \frac{2^{\frac{1}{\al}+1}}{\pi}
        \int_0^\infty
        \frac{u^{\frac{1}{\al}-1}}{\left(1+u^2\right)^{\frac{1}{\al}+1}}\,du
        =
        \frac{1}{\pi}\int_0^\pi \sin^{\frac{1}{\al}}y\,\cot(y/2)\,dy
        \in(0,\infty).
    \end{align*}
    Hence $h$ extends to a positive continuous function on $[0,\pi]$, as claimed.
\end{proof}

The full kernel satisfies the same bound.

\begin{corollary}\label{CorKernelTwoSidedFull}
    Let $\al>0$ and $n\geq3$. Then, there exist constants $0<c\leq C$, depending only on $\al$, such that
    \begin{align*}
        c\,\sin x\leq
        \int_0^\pi\sin^{\frac{1}{\al}}y\,\kernel(x,y)\,dy
        \leq C\sin x
        \qquad\text{for every }x\in[0,\pi].
    \end{align*}
\end{corollary}
\begin{proof}
    By Lemma~\ref{Lemma_Kernel_Def} and the definition of $\psi_n$ in Lemma~\ref{LemmaPhiIncreasing},
    \begin{align}\label{IdentityKernelSplit}
        \kernel(x,y)=
        \frac{1}{2\pi}\log\!\left(\frac{1-\cos(x+y)}{1-\cos(x-y)}\right)
        +\frac{1}{2\pi}\big(\psi_n(x+y)-\psi_n(x-y)\big),
    \end{align}
    From Lemma~\ref{LemmaEvenMonotone} and \eqref{IdentityKernelSplit}, we get
    \begin{align}\label{IneqKernelLowerBoundFromSplit}
        \kernel(x,y)\geq\frac{1}{2\pi}\log\!\left(\frac{1-\cos(x+y)}{1-\cos(x-y)}\right)
        \qquad\text{for every }(x,y)\in[0,\pi]^2.
    \end{align}
    The lower bound in the statement is therefore the one already proved in Lemma~\ref{LemmaKernelTwoSided}.

    For the upper bound,
    \begin{align}
        \frac{6}{\pi}\sum_{k\geq 1}\frac{\sin(kx)}{n^2k^3-4k}
        \int_0^\pi\sin^{\frac{1}{\al}}y\,\sin(ky)\,dy
        &\leq
        \frac{6}{\pi}\sum_{k\geq 1}\frac{|\sin(kx)|}{n^2k^3-4k}
        \int_0^\pi\sin^{\frac{1}{\al}}y\,dy\nonumber\\
        &\leq
        \frac{6}{\pi}\sum_{k\geq 1}\frac{1}{n^2k^2-4}
        \int_0^\pi\sin^{\frac{1}{\al}}y\,dy\,\sin x,\label{IneqUpperBoundSecondTerm}
    \end{align}
    where we have used $|\sin(ky)|\leq1$ and $|\sin(kx)|\leq k\sin x$, since
    \begin{align*}
        |\sin((k+1)x)|=|\sin(kx)\cos x+\cos(kx)\sin x|\leq|\sin(kx)|+\sin x,
    \end{align*}
    which gives $|\sin(kx)|\leq k\sin x$ by induction on $k$.

    From \eqref{IneqUpperBoundSecondTerm} and Lemma~\ref{LemmaKernelTwoSided}, we finally have
    \begin{align*}
        \int_0^\pi\sin^{\frac{1}{\al}}y\,\kernel(x,y)\,dy
        \lesssim
        \sin x
        +\frac{6}{\pi}\sum_{k\geq 1}\frac{1}{n^2k^2-4}
        \int_0^\pi\sin^{\frac{1}{\al}}y\,dy\,\sin x
        \lesssim
        \sin x.
    \end{align*}
\end{proof}

\subsection{Extending the uniqueness class}\label{extending}

Finally, we check that every nonnegative solution of \eqref{EqnIntegral} other than $0$ belongs to $\Kcal$. For $\al>1$ this upgrades the uniqueness in $\Kcal$ given by Proposition~\ref{PropGeneralUniqueness} to uniqueness among all nonnegative solutions.

\begin{lemma}\label{LemmaSolutionComparable}
    Let $\al>1$, let $n\geq3$ and let $v\in C([0,\pi])$ be nonnegative, not identically zero and such that $v$ satisfies \eqref{EqnIntegral}. Then $v\in\Kcal$.
\end{lemma}
\begin{proof}
    By Lemma~\ref{LemmaRegularityHolder}, $v\in C^{1}([0,\pi])$, and by Lemma~\ref{LemmaKernelPositive}, $v(0)=v(\pi)=0$ and $v>0$ on $(0,\pi)$. The mean value theorem then gives $v(x)\leq\norm{\dx v}_{C([0,\pi])}\min\{x,\pi-x\}$, and since $\min\{x,\pi-x\}\leq\tfrac{\pi}{2}\sin x$ on $[0,\pi]$, we obtain
    \begin{align*}
        v(x)\leq\frac{\pi}{2}\norm{\dx v}_{C([0,\pi])}\sin x\qquad\text{for every }x\in[0,\pi].
    \end{align*}

    For the lower bound, fix $[a,b]\subset(0,\pi)$ and write $s=\tan\big(\frac{x}{2}\big)$, $u=\tan\big(\frac{y}{2}\big)$. By \eqref{IneqKernelLowerBoundFromSplit} and \eqref{IdentityTangentKernel},
    \begin{align}\label{IneqKernelLowerBoundCompact}
        \kernel(x,y)\geq
        \frac{1}{\pi}\log\left|\frac{s+u}{s-u}\right|
        \geq\frac{2}{\pi}\,\frac{\min\{s,u\}}{\max\{s,u\}},
    \end{align}
    where we used in the second inequality that
    \begin{align*}
        \log\!\left(\frac{1+z}{1-z}\right)=2\sum_{k\geq0}\frac{z^{2k+1}}{2k+1}\geq2z
        \qquad\text{for all }0<z<1.
    \end{align*}

    We have $\tfrac{s}{1+s^2}\leq\min\{s,\tfrac1s\}$ for every $s>0$, so that
    \begin{align}\label{IneqMinBoundSinX}
        \min\{s,\tfrac1s\}\geq\frac{s}{1+s^2}=\frac12\sin x.
    \end{align}
    For $y\in[a,b]$ the variable $u$ ranges in a compact subinterval of $(0,\infty)$, so the last quotient in \eqref{IneqKernelLowerBoundCompact} is bounded from below by a multiple of $\min\{s,\tfrac{1}{s}\}$. Using \eqref{IneqKernelLowerBoundCompact} and \eqref{IneqMinBoundSinX}, we get $\kernel(x,y)\geq c_{a,b}\sin x$ for every $x\in[0,\pi]$ and every $y\in[a,b]$, and therefore
    \begin{align*}
        v(x)=\int_0^\pi v(y)^{\frac{1}{\al}}\,\kernel(x,y)\,dy
        \geq c_{a,b}\Big(\int_a^b v(y)^{\frac{1}{\al}}\,dy\Big)\sin x,
    \end{align*}
    where the integral is positive because $v>0$ on $(0,\pi)$.
\end{proof}

The following lemma will let us extend the uniqueness in $\Kcal$ to a wider family.

\begin{lemma}\label{LemmaRigidityL2}
    Let $\al>1$, let $n\geq3$, and let $f\in L^2(\T)$ be an odd $\frac{2\pi}{n}$-periodic solution of \eqref{EqnStatFolded}, in the sense of Remark~\ref{RemWeakSense}, such that $\OpSUnfold f$ does not vanish on $\big(0,\tfrac{\pi}{n}\big)$. Then there exists $\Gamma\neq0$ such that, setting $g(x):=\Gamma f\big(\tfrac{x}{n}\big)$ for $x\in[0,\pi]$, the function $v:=(-g)^\al$ satisfies \eqref{EqnIntegral} and $v\in\Kcal$.
\end{lemma}
\begin{proof}
    Let $f$ be as in the statement, and set $g_0(x):=f\big(\tfrac{x}{n}\big)$. We first show that $g_0\in H^1_{\mathrm{loc}}((0,\pi))$.

    In particular $f\neq 0$, and $g_0$ solves \eqref{EqnStatFoldedg} in the sense of Remark~\ref{RemWeakSense}. Moreover, by Lemma~\ref{Lemma_FourierUnfoldFold}, we have $\OpSUnfold f,\OpS g_0\in C(\T)$, and $\OpS g_0(x)=n\,\OpSUnfold f\big(\tfrac{x}{n}\big)$, so $\OpS g_0$ does not vanish on $(0,\pi)$.

    Fix a compact interval $J\subset(0,\pi)$ and test the weak form of \eqref{EqnStatFoldedg} (Remark~\ref{RemWeakSense}) on $J$ (legitimate by density),
    \begin{align*}
        -\al\int_\T g_0\,\OpS g_0\,\dx\!\left(\frac{\phi}{\OpS g_0}\right)dx
        =(\al+1)\int_\T g_0\,\dx(\OpS g_0)\,\frac{\phi}{\OpS g_0}\,dx
        \qquad\text{for every }\phi\in C_c^\infty(J).
    \end{align*}
    Applying the Leibniz rule, we find
    \begin{align*}
        -\al\int_\T g_0\,\dx\phi\,dx
        =\int_\T\frac{g_0\,\dx(\OpS g_0)}{\OpS g_0}\,\phi\,dx
        \qquad\text{for every }\phi\in C_c^\infty(J).
    \end{align*}
    Since the right hand side is integrable, we get
    \begin{align*}
        g_0\in W^{1,1}(J)\subset L^\infty(J),
    \end{align*}
    by the one-dimensional Sobolev embedding. Feeding this back into the right hand side,
    \begin{align*}
        \left\|\frac{g_0\,\dx(\OpS g_0)}{\OpS g_0}\right\|_{L^2(J)}
        \leq
        \left\|\frac{g_0}{\OpS g_0}\right\|_{L^\infty(J)}
        \|\dx(\OpS g_0)\|_{L^2(J)}
        <\infty,
    \end{align*}
    so the right hand side lies in $L^2(J)$, so we have the same for $\dx g_0$, and consequently $g_0\in H^1(J)$.
    Since $J\subset(0,\pi)$ was an arbitrary compact interval, this gives $g_0\in H^1_{\mathrm{loc}}((0,\pi))$, as claimed, which in turn gives $g_0\in C((0,\pi))$ by the one-dimensional Sobolev embedding again. Moreover, \eqref{EqnStatFoldedg} holds pointwise a.e. on $(0,\pi)$.

    We have $\dx|g_0|^\al=\al|g_0|^{\al-1}\sgn(g_0)\,\dx g_0$, so \eqref{EqnStatFoldedg} gives
    \begin{align*}
        \dx|g_0|^\al=\frac{|g_0|^\al}{\OpS g_0}\,\dx\OpS g_0
        \qquad\text{on }(0,\pi).
    \end{align*}
    Hence $|g_0|^\al/\OpS g_0$ has vanishing derivative there, so $|g_0|^\al=C\,\OpS g_0$ on $(0,\pi)$ for some constant $C$, and $C\neq 0$ because $f\neq 0$.

    Set $g:=\Gamma g_0$ with $\Gamma:=\sgn(C)|C|^{\frac{1}{1-\al}}$, well defined because $\al\neq 1$, so that $|g|^\al=\OpS g$. Since $\OpS g=\Gamma\,\OpS g_0$ is continuous on $\T$ and does not vanish on $(0,\pi)$, we have $|g|^\al=\OpS g\neq0$ there, so $g$ is nonvanishing on $(0,\pi)$ and, by continuity, has a constant sign. It must be negative, since $g>0$ would give $\OpS g<0$ by the positivity of $-\OpS$ from Lemma~\ref{LemmaKernelPositive}, contradicting $\OpS g=|g|^\al>0$. Hence $v:=(-g)^\al=\OpS g$ is continuous on $\T$, and $g=-v^{\frac{1}{\al}}$ on $(0,\pi)$. Since $v=\OpS g$, the continuity of $v$ extends to $[0,\pi]$, and we have that $v$ satisfies \eqref{EqnIntegral}. Finally, Lemma~\ref{LemmaSolutionComparable} gives $v\in\Kcal$.
\end{proof}

\section{Main theorem}

The ranges $\al>1$ and $\al\in\big(\tfrac{1}{2},1\big)$ are treated by different arguments and lead to different conclusions, so we state the corresponding results separately. For $\al=1$ the equation is trivially solvable, since any multiple of $\sin(nx)$ solves \eqref{EqnStatFolded}.

\begin{theorem}\label{TheoremAlBig}
    For any $\al>1$ and any $n\geq 3$, there exists an odd $\frac{2\pi}{n}$-periodic function $f_\al$, belonging to $C^{\frac{1}{\al}}(\T)$, such that its $n$-th Fourier coefficient does not vanish, i.e. $(\ha{f_\al})_n\neq 0$, and $f_\al$ solves equation \eqref{EqnStatFolded},
    \begin{align*}
        \al \OpSUnfold f_\al\, \dx f_\al=f_\al\, \dx\OpSUnfold f_\al.
    \end{align*}
    It is given by $f_\al(x)=\ti g_\al(nx)$, where $g_\al$ is as in Lemma~\ref{LemmaRegularityHolder}, $\ti g_\al$ is its odd extension to $\T$, and $v_\al$ is the only nonnegative solution of \eqref{EqnIntegral} other than $0$. Moreover, if $f\in L^2(\T)$ is an odd $\frac{2\pi}{n}$-periodic solution of \eqref{EqnStatFolded}, and $\OpSUnfold f$ does not vanish on $\big(0,\tfrac{\pi}{n}\big)$, then $f$ is a multiple of $f_\al$.
\end{theorem}
\begin{proof}
    The operator $-\OpS$ is linear and maps nonnegative functions to nonnegative functions by Lemma~\ref{LemmaKernelPositive}, and it satisfies hypothesis \eqref{HypTwoSided} by Corollary~\ref{CorKernelTwoSidedFull}. Proposition~\ref{PropGeneralUniqueness} therefore gives a unique fixed point $v_\al\in\Kcal$ of $v\mapsto-\OpS\big(v^{\frac{1}{\al}}\big)$. By Lemma~\ref{LemmaSolutionComparable}, every nonnegative solution of \eqref{EqnIntegral} other than $0$ belongs to $\Kcal$, so $v_\al$ is the only one.

    It remains to pass from $v_\al$ to $f_\al$. Since $v_\al\in\Kcal$, it is positive on $(0,\pi)$ and vanishes at $0$ and $\pi$, so Corollary~\ref{CorRegularitySolution} applies and gives every remaining assertion, namely that $f_\al$ is odd and $\frac{2\pi}{n}$-periodic, belongs to $C^{\frac{1}{\al}}(\T)$, solves \eqref{EqnStatFolded}, and satisfies $(\ha{f_\al})_n\neq 0$.

    For the last statement, let $f\in L^2(\T)$ be an odd $\frac{2\pi}{n}$-periodic solution of \eqref{EqnStatFolded} such that $\OpSUnfold f$ does not vanish on $\big(0,\tfrac{\pi}{n}\big)$. By Lemma~\ref{LemmaRigidityL2}, there exists $\Gamma\neq0$ such that, setting $g(x):=\Gamma f\big(\tfrac{x}{n}\big)$, the function $v:=(-g)^\al$ satisfies \eqref{EqnIntegral} and $v\in\Kcal$. The uniqueness established in the first paragraph gives $v=v_\al$, so $g=g_\al$ and $\Gamma f=f_\al$.
\end{proof}

\begin{theorem}\label{TheoremAlSmall}
    For any $\al\in\big(\tfrac{1}{2},1\big)$ and any $n\geq 3$, there exists an odd $\frac{2\pi}{n}$-periodic function $f_\al\in C^{1,\frac{1}{\al}-1}(\T)$ such that its $n$-th Fourier coefficient does not vanish, i.e. $(\ha{f_\al})_n\neq 0$, and $f_\al$ solves equation \eqref{EqnStatFolded}.
\end{theorem}
\begin{proof}
    By Proposition~\ref{PropGeneralExistence}, applied with $\OpSGen=-\OpS$ and with hypotheses checked in Section~\ref{checking}, there exists a nonzero fixed point $v_\al$ of $v\mapsto-\OpS\big(v^{\frac{1}{\al}}\big)$ that is a multiple of an element of $\ConvSet$.

    By Lemma~\ref{LemmaKernelPositive}, $v_\al$ is positive on $(0,\pi)$ and vanishes at $0$ and $\pi$, so Corollary~\ref{CorRegularitySolution} applies.
\end{proof}

\begin{proofsth}{Theorem~\ref{MainTheorem}}
    This follows from the particular case $\al=2$ of Theorem~\ref{TheoremAlBig}. Let $f:=f_2$ be the function given by the theorem. In particular, the profile of the stationary solution $\Theta$ is unique in the sense stated there. It remains to verify that $\Theta$ solves \eqref{SQG} in the sense of Remark~\ref{RemDistSense}.

    Write
    \begin{align*}
        u(x_1,x_2)=-\Lambda^{-1}\Theta(x_1,x_2)
        =r^2\OpSUnfold f(x).
    \end{align*}
    Since $\OpSUnfold=\dx^{-1}\!H$, we have $\OpSUnfold f\in C^{1,\frac{1}{2}}(\T)$, so $u\in C^1(\RR^2)$. Since $f$ is continuous, $\Theta\nabla^\perp u$ is then continuous on $\RR^2$, so it is locally integrable and the integral of Remark~\ref{RemDistSense} is absolutely convergent.

    In polar coordinates, for every $\phi(x_1,x_2)=\phi(r,x)$, in the radial-angular frame,
    \begin{align*}
        \nabla\phi=\Big(\pa_r\phi,\ \tfrac1r\pa_x\phi\Big)
        \qquad\text{and}\qquad
        \nabla^\perp\phi=\Big(-\tfrac1r\pa_x\phi,\ \pa_r\phi\Big).
    \end{align*}
    Hence,
    \begin{align*}
        \int_{\RR^2}\Big(\Theta\,\nabla^\perp u\cdot\nabla\psi\Big)(x_1,x_2)\,dx
        &=
        \int_0^\infty\intpi
        2r^2 f(x)\,\OpSUnfold f(x)\,\pa_x\psi(r,x)
        \,dx\,dr\\
        &\quad-
        \int_0^\infty\intpi
        r^3 f(x)\,\pa_x\OpSUnfold f(x)\,\pa_r\psi(r,x)
        \,dx\,dr.
    \end{align*}
    Integrating by parts the first term in $x$ and the second term in $r$, we obtain
    \begin{align}
        \int_{\RR^2}\Big(\Theta\,\nabla^\perp u\cdot\nabla\psi\Big)(x_1,x_2)\,dx
        &=
        \int_0^\infty\intpi
        r^2 f(x)\,\pa_x\OpSUnfold f(x)\,\psi(r,x)
        \,dx\,dr\nonumber\\
        &\quad-
        \int_0^\infty\intpi
        2r^2\OpSUnfold f(x)\,\pa_xf(x)\,\psi(r,x)
        \,dx\,dr,\label{EqnMainThmPairing}
    \end{align}
    with no boundary contribution in either integration by parts, in $x$ because $f\,\OpSUnfold f\,\psi$ is continuous on $\T$, and in $r$ because $\psi$ has compact support and the integrand vanishes at $r=0$. Finally, by \eqref{EqnStatFolded} with $\al=2$, understood in the sense of Remark~\ref{RemWeakSense}, the right-hand side of \eqref{EqnMainThmPairing} vanishes.
\end{proofsth}

For a $2\pi$-periodic function $f$, the stationary version of \eqref{DeGregGen} reads as
\begin{equation}\label{EqnStatDeGreg}
    \al u\, \dx f = f\, \dx u,\qquad \dx u=Hf.
\end{equation}

\begin{theorem}\label{theorem:dgreggen}
    For any $\al>\tfrac{1}{2}$ and any $n\geq 1$, there exists an odd $\frac{2\pi}{n}$-periodic function $f_\al$, belonging to $C^{\frac{1}{\al}}(\T)$ if $\al>1$ and to $C^{1,\frac{1}{\al}-1}(\T)$ if $\al<1$, such that its $n$-th Fourier coefficient does not vanish, i.e. $(\ha{f_\al})_n\neq 0$, and $f_\al$ solves the stationary generalized De Gregorio equation \eqref{EqnStatDeGreg}. For $\al>1$, the function $f_\al$ arises from a fixed point that is unique in the sense of Theorem~\ref{TheoremAlBig}.
\end{theorem}

\begin{proof}
    It suffices to treat $n=1$, since the case $n\geq 2$ follows by constructing the solution for $n=1$ and rescaling by $x\mapsto nx$, as in Remark~\ref{RemRescaling}. For $\al=1$, any multiple of $\sin(nx)$ is a solution.

    For $\al\neq 1$, the proof follows exactly the proofs of Theorems~\ref{TheoremAlBig} and \ref{TheoremAlSmall}, substituting $\OpS$ by $-\dx^{-1}\!H=\Lambda^{-1}$, whose kernel is the logarithmic term of $\kernel$ alone. It is linear and maps nonnegative functions to nonnegative functions by the proof of Lemma~\ref{LemmaKernelPositive}, with $\phi(x)=\log(1-\cos x)$ in place of $\phi_n$, for which $\dx\phi(x)=\cot\big(\tfrac{x}{2}\big)>0$ on $(0,\pi)$ directly gives $\phi$ increasing, so Lemma~\ref{LemmaPhiIncreasing} is not needed and Lemma~\ref{LemmaEvenMonotone} applies directly.

    For $\al>1$, hypothesis \eqref{HypTwoSided} of Proposition~\ref{PropGeneralUniqueness} is Lemma~\ref{LemmaKernelTwoSided} itself, with no Fourier part left to estimate. Proposition~\ref{PropGeneralUniqueness} therefore gives a fixed point $v_\al$, unique in the sense of Theorem~\ref{TheoremAlBig}.

    For $\al\in\big(\tfrac{1}{2},1\big)$, hypothesis~\ref{enume_first} of Proposition~\ref{PropGeneralExistence} holds with $\rho=1$, since the multiplier of $\Lambda^{-1}$ equals $1$ at $k=1$. For hypothesis~\ref{enume_second}, the nonnegativity of the kernel is the one just proved, and the remaining inequality follows as in the proof of Lemma~\ref{LemmaBoundAlpha}, where the $\ep$-$\de$ condition of Lemma~\ref{LemmaEpDeltaGood} is now obtained in the same way as in the proof of that lemma, without having to bound the Fourier (compact) part, which there only contributes an $O(\de)$ term. Proposition~\ref{PropGeneralExistence} then gives a nonzero fixed point $v_\al$ of $v\mapsto\Lambda^{-1}\big(v^{\frac{1}{\al}}\big)$ that is a multiple of an element of $\ConvSet$.

    The argument of Corollary~\ref{CorRegularitySolution} applies with $\OpS$ replaced by $\Lambda^{-1}$, since the operator has no smoothing part left to bound, so in either case $v_\al$ gives a function $f_\al$ that solves \eqref{EqnStatDeGreg} and has the stated regularity.
\end{proof}

\appendix
\section{Auxiliary lemmas}
This appendix contains the proofs of the auxiliary lemmas used in the proof of Lemma~\ref{LemmaBoundAlpha}, as well as a further intermediate lemma used in their proofs.

\begin{proofsth}{Lemma~\ref{LemmaKernelBoundSine}} The kernel $\kernel$ has a logarithmic part and a Fourier series part (see Lemma~\ref{Lemma_Kernel_Def}). To bound the latter one, we apply the symmetry condition \eqref{symWiKernel} to $\ti{K}_n(u,z)=\kernel\big(\tfrac{u+z}{2},\tfrac{u-z}{2}\big)$,
    \begin{align*}
        \kernel(\pi-x,\pi-y)=
        \ti{K}_n(2\pi-x-y,y-x)=
        \ti{K}_n(x+y,x-y)=\kernel(x,y).
    \end{align*}
    Hence, it suffices to prove the bound for $y \in [0,\tfrac{\pi}{2}]$, since the case $y \in [\tfrac{\pi}{2},\pi]$ follows by symmetry.

    For $y \in [0,\tfrac{\pi}{2}]$, we have $y\leq\pi\sin y$. Consequently, the Fourier series part of $\kernel$ can be estimated as
    \begin{align}
    \label{IneqFourierPart}
        \sum_{k\geq 1}\frac{1}{n^2k^3-4k} \sin(kx)\sin(ky) =
        \int_0^y\sum_{k\geq 1}\frac{1}{n^2k^2-4} \sin(kx)\cos(ks) \,ds
        \lesssim
        \int_0^y \,ds
        \leq
        y\leq \pi\sin y.
    \end{align}

    It remains to control the logarithmic term. We claim that
    \begin{align}
    \label{ineqn_lemma_explcit}
        \log\!\left(\frac{1-\cos(x+y)}{1-\cos(x-y)}\right)\leq
        2\sin y\,\cot\!\left(\frac{|x-y|}{2}\right).
    \end{align}

    The kernel $\kernel(x,y)$ is not defined at $x=y$, as noted in the proof of Lemma~\ref{LemmaKernelPositive}. For $x\neq y$, $\te := x-y$ satisfies $\te \in [-\pi,\pi]-\{0\}$, so $1-\cos\te>0$. Then, inequality \eqref{ineqn_lemma_explcit} becomes
    \begin{align}
    \label{IneqLog}
        \log\!\left(\frac{1-\cos\te\,\cos(2y)+\sin\te\,\sin(2y)}{1-\cos\te}\right)\leq
        2\sin y\,\frac{\sin|\te|}{1-\cos\te}=:\Gamma.
    \end{align}
    Since $\Gamma\geq 0$, we have
    \begin{align*}
        1+\Gamma+\frac{1}{2}\Gamma^2 \leq \exp(\Gamma).
    \end{align*}
    Therefore, to prove inequality \eqref{IneqLog}, it suffices to check that
    \begin{align}\label{QuadraticIneq}
        \frac{1-\cos\te\,\cos(2y)+\sin\te\,\sin(2y)}{1-\cos\te}
        \leq
        1+\Gamma+\frac{1}{2}\Gamma^2.
    \end{align}

    We now simplify \eqref{QuadraticIneq} through a chain of equivalent inequalities. Multiplying both sides by $1-\cos\te>0$, we obtain
    \begin{align*}
        1-\cos\te\,\cos(2y)+\sin\te\,\sin(2y)
        \leq
        1-\cos\te+
        2\sin y\,\sin|\te| + 2\sin^2 y
        \frac{1-\cos^2\te}{1-\cos \te}.
    \end{align*}
    Now, we subtract $1$ and use the identity $\frac{1-\cos^2\te}{1-\cos\te}=1+\cos\te$, so \eqref{QuadraticIneq} is equivalent to
    \begin{align}\label{EquiQuad1}
        -\cos\te\,\cos(2y)
        +
        \sin\te\, \sin(2y)
        \leq
        -\cos\te+
        2\sin y\,\sin|\te| + 2\sin^2 y\,(1+\cos\te).
    \end{align}
    Observe that
    \begin{align}\label{Trigoidenty}
        -\cos\te\,\cos(2y)=-\cos\te+2\sin^2 y\,\cos\te.
    \end{align}

    Substituting \eqref{Trigoidenty} into \eqref{EquiQuad1}, we obtain
    \begin{align}\label{EquiQuad2}
        \sin\te\, \sin(2y)
        \leq
        2\sin y\,\sin|\te| + 2\sin^2 y.
    \end{align}
    Inequality \eqref{EquiQuad2} now follows from
    \begin{align*}
        \sin\te\,\cos y\leq
        \sin|\te|+\sin y,
    \end{align*}
    after multiplying by $2\sin y\geq 0$. This last inequality holds because $\sin\te\,\cos y\leq|\sin\te|=\sin|\te|$ for $\te\in[-\pi,\pi]$, and $\sin y\geq 0$ for $y\in[0,\pi]$. The statement then follows from \eqref{IneqFourierPart} and \eqref{ineqn_lemma_explcit}.
\end{proofsth}

\begin{proofsth}{Lemma~\ref{Lemma_when_isnegative}}
    The derivative of $\beta$ is given by
    \begin{align*}
        \pa_M\beta(M) = \al^{-1} A_1 M^{\frac{1}{\al}-1}-1.
    \end{align*}
    Hence, $\beta$ has a unique critical point $M_0$, determined by $M_0^{1-\frac{1}{\al}}=A_1\al^{-1}$.

    The second derivative satisfies
    \begin{align*}
        \pa^2_M\beta(M) = A_1\al^{-1} (\al^{-1}-1)M^{\frac{1}{\al}-2}>0.
    \end{align*}
    Hence, $M_0$ is a minimum of $\beta$.

    Evaluating $\beta$ at $M_0$, we obtain
    \begin{align*}
        \beta(M_0)=
        M_0\big(A_1 M_0^{\frac{1}{\al}-1}-1\big) + A_2=M_0(\al-1) +A_2.
    \end{align*}

    Hence, the minimum is negative if and only if $A_2<(1-\al)M_0$. Substituting $M_0=(A_1\al^{-1})^{\frac{1}{1-\al^{-1}}}$ and raising both sides to the power $1-\al>0$ and multiplying by $A_1^\al$, this is equivalent to
    \begin{align*}
        A_1^\al A_2^{1-\al} < (1-\al)^{1-\al} \al^\al.
    \end{align*}
\end{proofsth}

\begin{lemma}\label{LemmaClauEstimate}
    Given any $\de\in(0,\pi)$, let $\clau$ be the Clausen function introduced in Lemma~\ref{LemmaSineSeriesPositive} and let $\rho_\de:[0,\pi]\to\RR$ be the function given by
    \begin{align*}
        \rho_\de(x):=
        \clau(2x-\de) - \clau(2x+\de).
    \end{align*}
    Then, $\rho_\de$ is increasing on $[0,\tfrac{\pi}{2}]$ and attains its maximum at $\tfrac{\pi}{2}$.
\end{lemma}
\begin{proof}
    First, observe that $\rho_\de$ is symmetric with respect to $\tfrac{\pi}{2}$. Indeed,
    \begin{align*}
        \rho_\de(\pi-x) &=
        \clau(2(\pi-x)-\de) - \clau(2(\pi-x)+\de)\\
        &=
        \clau(-2x-\de) - \clau(-2x+\de)=
        \clau(2x-\de) - \clau(2x+\de)
        =\rho_\de(x),
    \end{align*}
    where we have used that $\clau$ is odd and $2\pi$-periodic. Hence, it is enough to study $\rho_\de$ on $[0,\tfrac{\pi}{2}]$.

    For $x\neq\tfrac{\de}{2}$, neither $2x-\de$ nor $2x+\de$ is an integer multiple of $2\pi$. Differentiating $\rho_\de$ and using \eqref{clauder}, we obtain
    \begin{align*}
        \dx\rho_\de(x) =
        2\dx\clau(2x-\de) - 2\dx\clau(2x+\de)=
        \log\!\left(\frac{1-\cos(2x+\de)}{1-\cos(2x-\de)}\right).
    \end{align*}
    Thus, $\dx\rho_\de(x)\geq 0$ whenever
    \begin{align}\label{CosIneq}
        \cos(2x-\de)\geq\cos(2x+\de).
    \end{align}

    To check \eqref{CosIneq}, note that
    \begin{align*}
        \cos(2x-\de)-\cos(2x+\de)
        =
        2\sin(2x)\sin\de.
    \end{align*}
    Since $x\in[0,\tfrac{\pi}{2}]$, we have $2x\in[0,\pi]$ and hence $\sin(2x)\geq 0$, while $\de\in(0,\pi)$ gives $\sin\de>0$. The difference is therefore nonnegative, so $\dx \rho_\de(x)\geq 0$ for every $x\in[0,\tfrac{\pi}{2}]$ with $x\neq\tfrac{\de}{2}$. Since $\rho_\de$ is continuous, it is increasing on this interval, and by symmetry it attains its maximum at $x=\tfrac{\pi}{2}$.
\end{proof}

We prove Lemma~\ref{LemmaEpDeltaGood} before Lemma~\ref{LemmaLimitAlSmallIsZero}, since the latter's proof relies on it.

\begin{proofsth}{Lemma~\ref{LemmaEpDeltaGood}}
    Let $I_{\de,x}$ be the interval $[0,\pi]\cap(x-\de,x+\de)$, so that
    \begin{align*}
        \Ical(\de)=\sup_{x\in[0,\pi]}\int_{I_{\de,x}}
        \kernel(x,y) \,dy.
    \end{align*}

    Recall the expression for $\kernel$ obtained in Lemma~\ref{Lemma_Kernel_Def},
    \begin{align*}
        \kernel(x,y)=
        \frac{1}{2\pi}\log\!\left(\frac{1-\cos(x+y)}{1-\cos(x-y)}\right) +
        \frac{6}{\pi}\sum_{k\geq 1}
        \frac{1}{n^2k^3-4k}\sin(kx)\sin(ky).
    \end{align*}
    We bound the logarithmic and Fourier parts of the kernel $\kernel$ separately.

    From \eqref{equ_expan_log} we obtain
    \begin{align}
        \log\!\left(\frac{1-\cos(x+y)}{1-\cos(x-y)}\right)=
        2\sum_{k\geq 1} \frac{1}{k}\Big(\cos(k(x-y))-\cos(k(x+y))\Big)=
        4\sum_{k\geq 1}\frac{1}{k}\sin(kx)\sin(ky).
        \label{log_to_fou}
    \end{align}
    Accordingly, we define
    \begin{align*}
        F(x):=
        \int_{I_{\de,x}}
        \sum_{k\geq 1}\frac{1}{k}
        \sin(kx)\sin(ky)
        \,dy.
    \end{align*}
    We claim that $F$ attains its maximum at $x=\tfrac{\pi}{2}$, which we prove in four steps.

    First, $F$ is symmetric with respect to $\tfrac{\pi}{2}$. Indeed, we have
    \begin{align*}
        \sin\big(k(\pi-x)\big)\sin\big(k(\pi-y)\big)=
        (-1)^{2k+2}\sin(kx)\sin(ky)=
        \sin(kx)\sin(ky).
    \end{align*}
    Moreover, the change of variables $y\mapsto\pi-y$ maps $I_{\de,\pi-x}$ onto $I_{\de,x}$. Therefore, $F(\pi-x)=F(x)$, and it is enough to study $F$ on $[0,\tfrac{\pi}{2}]$.

    Second, $F$ is increasing on $[\de,\tfrac{\pi}{2}]$. For $x\in[\de,\tfrac{\pi}{2}]$, we have $I_{\de,x}=(x-\de,x+\de)$. Integrating $F$ term by term,
    \begin{align}
        F(x)&=
        \sum_{k\geq 1}\frac{1}{k^2}
        \sin(kx)[\cos(k(x-\de)) - \cos(k(x+\de))]
        \nonumber\\&=
        \sum_{k\geq 1}\frac{1}{k^2}\left[\frac{1}{2}\sin(k(2x-\de))-\frac{1}{2}\sin(k(2x+\de))+\sin(k\de)\right]
        \nonumber\\&=
        \clau(\de)+\frac{1}{2}\clau(2x-\de)-\frac{1}{2}\clau(2x+\de)=
        \clau(\de)+\frac{1}{2}\rho_{\de}(x),
        \label{EqnLogAndrho}
    \end{align}
    where $\rho_\de$ and $\clau$ are the functions defined in Lemma~\ref{LemmaClauEstimate}. By that lemma, $\rho_\de$ is increasing on $[0,\tfrac{\pi}{2}]$, and hence so is $F$.

    Third, $F$ is increasing on $[0,\de]$. For $x\in[0,\de]$, we have $I_{\de,x}=[0,x+\de)$ and
    \begin{align}
        F(x)&=
        \sum_{k\geq 1}\frac{1}{k^2}
        \sin(kx)\big[1-\cos(k(x+\de))\big]
        \nonumber\\&=
        \clau(x)-\frac{1}{2}\clau(2x+\de)+\frac{1}{2}\clau(\de),
        \quad x\in[0,\de].
        \label{EqnFleft}
    \end{align}
    Differentiating \eqref{EqnFleft} with respect to $x$ and using \eqref{clauder}, we obtain
    \begin{align*}
        \dx F(x)=
        \frac{1}{2}\log\!\left(\frac{2-2\cos(2x+\de)}{2-2\cos x}\right).
    \end{align*}
    For $x\in(0,\de]$, we have $0< x\leq 2x+\de\leq 3\de\leq\pi$, and $\cos$ is decreasing on $[0,\pi]$, so the argument of the logarithm is at least $1$ and $\dx F(x)\geq 0$. Here, we have used that $\de<\tfrac{\pi}{3}$. Since $F$ is continuous on $[0,\de]$, it is increasing there.

    Fourth, both expressions, \eqref{EqnLogAndrho} and \eqref{EqnFleft}, coincide at $x=\de$. Hence, $F$ is continuous on $[0,\tfrac{\pi}{2}]$ and increasing there. Together with the symmetry from the first step and using $\clau(\pi+\de)=-\clau(\pi-\de)$, we find that
    \begin{align}
        \sup_{x\in[0,\pi]}F(x)=
        F\!\left(\tfrac{\pi}{2}\right)=
        \clau(\de)+\frac{1}{2}\rho_\de\!\left(\tfrac{\pi}{2}\right)=
        \clau(\de)+\clau(\pi-\de).
        \label{EqnLogAndClau}
    \end{align}

    Observe now that
    \begin{align}
        \clau(x)=
        \sum_{k\geq 1}\frac{1}{k^2}\sin(kx)&=
        \sum_{k\geq 1}\frac{1}{(2k-1)^2}\sin((2k-1)x)+
        \sum_{k\geq 1}\frac{1}{(2k)^2}\sin(2kx)\nonumber
        \\&=
        \sum_{k\geq 1}\frac{1}{(2k-1)^2}\sin((2k-1)x)+
        \frac{1}{4}\clau(2x).
        \label{EqnClauOddEven}
    \end{align}
    Consequently,
    \begin{align}
        \clau(\de)+\clau(\pi-\de)
        &=
        \sum_{k\geq 1}\frac{1}{k^2}\sin(k\de)+
        \sum_{k\geq 1}\frac{1}{k^2}\sin\big(k(\pi-\de)\big)\nonumber \\
        &=
        2\sum_{k\geq 1}\frac{1}{(2k-1)^2} \sin\big((2k-1)\de\big)=
        2 \clau(\de) - \frac12 \clau(2\de)
        \label{eqnrewriclaudiff},
    \end{align}
    where we have used \eqref{EqnClauOddEven} in the last equality. Therefore, from \eqref{clauder} and \eqref{eqnrewriclaudiff}, we have that
    \begin{align*}
        \pa_\de[\clau(\de)+\clau(\pi-\de)]
        &=
        \pa_\de\!\left[2\clau(\de)-\frac12\clau(2\de)\right]
        =
        \frac{1}{2}
        \log\!\left(\frac
        {2-2\cos(2\de)}
        {(2-2\cos\de)^2}\right)
        =
        \log\cot(\de/2).
    \end{align*}
    Since $\cos(\de/2)\leq 1$ and $\sin(\de/2)\geq\de/\pi$ for every $\de\in(0,\pi)$, we have
    \begin{align}\label{AsympClauDiff}
        \pa_\de[\clau(\de)+\clau(\pi-\de)]\leq
        -\log\de + \log\pi.
    \end{align}

    Hence, using \eqref{log_to_fou} and combining estimates \eqref{EqnLogAndClau} and \eqref{AsympClauDiff} with the identity $\clau(0)+\clau(\pi)=0$, we find that
    \begin{align}\label{ClauExplicit}
        \sup_{x\in[0,\pi]}\int_{I_{\de,x}}
        \log\!\left(\frac{1-\cos(x+y)}{1-\cos(x-y)}\right)\,dy\leq
        4\int_0^\de\log\cot(t/2)\,dt\leq
        -4\de\log\de +
        4\de(1+\log\pi).
    \end{align}

    Next, we estimate the Fourier series. For $n\geq 3$ and every $k\geq 1$, we have
    \begin{align*}
        n^2k^3-4k\geq 9k^3-4k=k\big(9k^2-4\big)\geq 5k>0,
    \end{align*}
    so the series is bounded uniformly in $x$ and $y$ by
    \begin{align*}
        \left|\sum_{k\geq 1}
        \frac{\sin(kx)\sin(ky)}{n^2k^3-4k}\right|\leq
        \sum_{k\geq 1}
        \frac{1}{9k^3-4k}<\infty.
    \end{align*}
    Since the length of $I_{\de,x}$ is at most $2\de$, this yields
    \begin{align}
        \sup_{x\in[0,\pi]}\int_{I_{\de,x}}
        \sum_{k\geq 1}
        \frac{\sin(kx)\sin(ky)}{n^2k^3-4k}
        \,dy\lesssim\de.
        \label{ClauNonExplicit}
    \end{align}

    Finally, from inequalities \eqref{ClauExplicit} and \eqref{ClauNonExplicit}, we deduce that
    \begin{align*}
        \Ical(\de) \lesssim
        \sup_{x\in[0,\pi]}\int_{I_{\de,x}}
        \log\!\left(\frac{1-\cos(x+y)}{1-\cos(x-y)}\right)\,dy
        +
        \sup_{x\in[0,\pi]}\int_{I_{\de,x}}
        \sum_{k\geq 1}
        \frac{\sin(kx)\sin(ky)}{n^2k^3-4k}
        \,dy\lesssim \de(1-\log\de).
    \end{align*}
\end{proofsth}

\begin{proofsth}{Lemma~\ref{LemmaLimitAlSmallIsZero}}
    We are computing a limit as $\de\to0^+$, so we can assume $\de<\tfrac{\pi}{3}$. Using the estimate in Lemma~\ref{LemmaEpDeltaGood}, we obtain
    \begin{align*}
    \lim_{\de\to0^+}
        \Ical(\de)^\al
        \left(1+\cot(\de/2)\right)^{1-\al}
        \lesssim
            \lim_{\de\to0^+}
            \frac{\left(\de(1-\log\de)\right)^\al}{\left(1+\cot(\de/2)\right)^{\al-1}}
        \lesssim
            \lim_{\de\to0^+}
            \left(
                \de(1-\log\de)
            \right)^\al \de^{\al-1}.
    \end{align*}
    We conclude that, for $\de>0$ small,
    \begin{align*}
        \Ical(\de)^\al
        \left(1+\cot(\de/2)\right)^{1-\al}\lesssim
        \de^{2\al-1}(-\log\de)^\al=
        \left(
        -\de^{2-\al^{-1}}\log\de\right)^\al.
    \end{align*}
    The right hand side tends to $0$ when $\de\to0^+$, because
    \begin{align*}
        \lim_{\de\to 0^+} \de^b \log\de = 0 \qquad \text{for all } b>0,
    \end{align*}
    and here $b=2-\al^{-1}$, which is positive if and only if $\al>\tfrac{1}{2}$. The left hand side is nonnegative, since $\kernel\geq 0$ by Lemma~\ref{LemmaKernelPositive} gives $\Ical\geq 0$, and $1+\cot(\de/2)>0$ for $\de$ small. The limit is therefore zero.
\end{proofsth}

\section*{Acknowledgments}
The author acknowledges financial support from the Severo Ochoa Programme for Centres of Excellence Grants \nolinkurl{CEX2019-000904-S} and \nolinkurl{CEX2023-001347-S}, and from grant \nolinkurl{PID2020-114703GB-I00}, funded by \nolinkurl{MCIN/AEI/10.13039/501100011033}. The author also thanks \'Angel Castro for his valuable comments and discussions during the preparation of this work.

\section*{Generative AI statement}
During the preparation of this work, the author used Claude (Anthropic) to assist with
copyediting and language revision, and to review the mathematical arguments and computations in the
manuscript for internal consistency and possible gaps. During this process, the ideas behind
the uniqueness argument, the refined regularity of the solution, and the proof of Lemma~\ref{LemmaPhiIncreasing} were worked out through discussion.
The author reviewed and verified the content of this publication and takes full responsibility for it.

\printbibliography

\begin{flushleft}
	\begin{minipage}{\linewidth}
	Miguel M.G. Pascual-Caballo\\
	\textsc{Instituto de Ciencias Matemáticas\\
		28049 Madrid, Spain}\\
	\textit{E-mail address:} miguel.martinez@icmat.es
	\end{minipage}
\end{flushleft}

\end{document}